\documentclass[10pt, reqno]{amsart}
\usepackage[utf8]{inputenc}

\usepackage{style}
\hypersetup{pdftitle={Ramsey numbers of loose triangles versus cliques}}

\newcommand{\ind}{\mathbf 1}
\newcommand{\dd}{\mathrm d}
\renewcommand{\phi}{\varphi}

\title{Ramsey numbers of loose triangles versus cliques}
\author{Ayush Basu}
\address{Department of Mathematics, University of Illinois Urbana--Champaign}
\email{ayushb@illinois.edu}
\author{Abhishek Dhawan}
\email{adhawan2@illinois.edu}
\thanks{Abhishek Dhawan is supported by NSF RTG grant DMS-1937241.}
\author{Abhishek Methuku}
\email{methuku@illinois.edu}
\thanks{Abhishek Methuku is supported by a Simons award SFI-MPS-TSM-00025583 and the UIUC Campus Research Board Award RB25050.}
\date{}
\begin{document}\raggedbottom

\begin{abstract}
An $r$-uniform loose triangle $C_3$ consists of three edges that intersect pairwise in exactly one vertex and have no vertex common to all three.
For graphs, this is the usual triangle, and the Ramsey number $R(C_3,K_t)$ has order of magnitude $t^2/\log t$ by celebrated results of Ajtai, Koml\'os, and Szemer\'edi and of Kim.

Kostochka, Mubayi, and Verstra\"ete proved that $R(C_3,K_t^{(r)})=O_r(t^{3/2})$ for every fixed $r\ge3$ and showed that this bound is sharp up to a logarithmic factor.
Motivated by the logarithmic saving in the graph case, they conjectured in 2013 that the upper bound could be improved to $o(t^{3/2})$.
This conjecture has since received significant attention and has been reiterated in several papers.
We prove the conjecture by showing that, for every fixed $r\ge3$,
\[
    R(C_3,K_t^{(r)})
    =O_r\left(\frac{t^{3/2}}{\sqrt{\log\log t}}\right).
\]
\end{abstract}

\maketitle\enlargethispage{3pt}

\section{Introduction}
\subsection{Ramsey numbers}\label{sec:background}
Quantitative aspects of Ramsey's theorem~\cite{Ramsey} were first studied by Erd\H{o}s and Szekeres~\cite{ES} for graphs and Erd\H{o}s and Rado~\cite{ER} for hypergraphs. A hypergraph is \emph{$r$-uniform}, or an \emph{$r$-graph}, if all of its edges have size $r$. For $r$-graphs $F$ and $G$, let $R(F,G)$ denote the smallest integer $N$ such that every red/blue coloring of the edges of $K_N^{(r)}$ contains a red copy of $F$ or a blue copy of $G$. A set of vertices is \emph{independent} if it contains no edge, and the \emph{independence number} $\alpha(H)$ is the largest size of an independent set in $H$. Then $R(F,K_t^{(r)})\le N$ if and only if every $F$-free $r$-graph on $N$ vertices has an independent set of size $t$. Thus, studying the growth of the \emph{off-diagonal Ramsey numbers} $R(F,K_t^{(r)})$ as $t\to\infty$ amounts to understanding how the structure of $F$ influences the size of independent sets in $F$-free hypergraphs.

For graphs, the classical bound of Erd\H{o}s and Szekeres~\cite{ES} gives $R(K_s,K_t)\le \binom{s+t-2}{s-1}=O(t^{s-1})$ for every fixed $s$, and hence $R(F,K_t)$ grows at most polynomially in $t$ for every fixed graph $F$. For $s=3$, Ajtai, Koml\'os, and Szemer\'edi~\cite{AKS} established the upper bound $R(K_3,K_t)=O(t^2/\log t)$, which Kim~\cite{Kim} famously showed to be tight up to a constant factor. The correct polynomial exponent for $R(K_4,K_t)$ remained unknown until the breakthrough of Mattheus and Verstra\"ete~\cite{MV}, who showed that $R(K_4,K_t)=\Omega(t^3/\log^4 t)$. Recently, Brada\v{c}~\cite{Bradac2026} determined the correct polynomial exponent for every fixed $s\ge 5$, showing that $R(K_s,K_t)=t^{s-1+o(1)}$.

Kostochka, Mubayi, and Verstra\"ete~\cite{KMV} initiated the study of an analogue of triangles for $r$-uniform hypergraphs motivated by the celebrated works of~\cite{AKS} and~\cite{Kim} on $R(K_3,K_t)$ for graphs. A natural analogue of a triangle in hypergraphs~\cite{MubayiRodl} is the \emph{loose triangle}, which consists of three distinct edges $e,f,g$ such that
\begin{equation}\label{eq: loose triangle def}
    |e\cap f|=|f\cap g|=|g\cap e|=1
    \qquad\text{and}\qquad
    e\cap f\cap g=\varnothing.
\end{equation}
For $r=2$, this is exactly $K_3$.
We denote the $r$-uniform loose triangle by $C_3$, with $r$ understood from context.
The loose triangle is the shortest loose cycle, and the work of Kostochka, Mubayi, and Verstra\"ete also inspired the study of Ramsey numbers of longer loose cycles and of Berge cycles versus cliques~\cite{Meroueh,NV}, as well as of tight cycles versus cliques~\cite{MubayiRodl}. For the loose triangle, they proved that, for every fixed $r\ge3$, there are constants $a_r,b_r,c_r>0$ such that, for all sufficiently large $t$,
\[
    \frac{a_r\,t^{3/2}}{(\log t)^{3/4 + b_r/\sqrt{\log \log t}}}
    \le R(C_3,K_t^{(r)})\le c_r\,t^{3/2},
\]
and for $r = 3$, they obtained the stronger lower bound $\Omega(t^{3/2}/(\log t)^{3/4})$.

In view of the logarithmic saving in the graph Ramsey number $R(K_3,K_t)=O(t^2/\log t)$ shown by~\cite{AKS}, Kostochka, Mubayi, and Verstra\"ete~\cite[Conjecture~1.1]{KMV} conjectured that the upper bound $O(t^{3/2})$ for loose triangles can also be improved by a factor that tends to infinity with $t$.
\begin{conjecture}[{Kostochka--Mubayi--Verstra\"ete~\cite{KMV}}]\label{conj: KMV}
    For every fixed integer $r\ge3$,
    \[
        R(C_3,K_t^{(r)})=o(t^{3/2})
        \qquad\text{as }t\to\infty.
    \]
\end{conjecture}
Conjecture~\ref{conj: KMV} has attracted considerable interest~\cite{MMNV,MubayiSuk}, motivating a line of work on Ramsey numbers and coloring problems for other hypergraph generalizations of the triangle~\cite{BMP,CM,LP,Meroueh}.

In particular, Cooper and Mubayi~\cite{CM} showed that forbidding $C_3$ together with the $3$-graphs $F_5=\{abc,abd,cde\}$ and $K_4^{(3)-}=\{abc,abd,acd\}$ (see also the discussion of Berge triangles below) yields an upper bound of $O(t^{3/2}/\sqrt{\log t})$. 
Our main result confirms Conjecture~\ref{conj: KMV} for every $r\ge3$, giving the first improvement on the upper bound $O(t^{3/2})$ for the loose triangle alone.

\begin{tcolorbox}[enhanced, breakable]
    \begin{theorem}\label{thm:Ramsey}
        For every integer $r\ge3$, there is a constant $C_r>0$ such that, for all $t\ge3$,
        \[
            R(C_3,K_t^{(r)})\le C_r\frac{t^{3/2}}{\sqrt{\log\log t}}.
        \]
    \end{theorem}
\end{tcolorbox}
The factor $\sqrt{\log\log t}$ arises naturally in our proof, as we explain in Section~\ref{sec:overview}.
The correct order of magnitude of $R(C_3,K_t^{(r)})$ remains unknown, and it would be interesting to determine whether one can match the lower bound construction of Kostochka, Mubayi, and Verstra\"ete~\cite{KMV}; for $r=3$, Mattheus, Mubayi, Nie, and Verstra\"ete~\cite{MMNV} suggested that this is the case, that is, that $R(C_3,K_t^{(3)})=\Theta\big(t^{3/2}/(\log t)^{3/4}\big)$.

Conjecture~\ref{conj: KMV} is part of a broader program, dating back to the work of Erd\H{o}s and Hajnal~\cite{ErdosHajnal72}, of understanding the growth of off-diagonal Ramsey numbers of hypergraphs. For $3$-uniform hypergraphs, they showed that $R(K_4^{(3)},K_t^{(3)})$ grows at least exponentially in $t$, a bound subsequently strengthened to $2^{\Omega(t\log t)}$ by Conlon, Fox, and Sudakov~\cite{CFS}. On the other hand, their work also implies that $R(F,K_t^{(3)})$ grows polynomially whenever $F$ is contained in an iterated blowup of an edge, in particular when $F$ is tripartite. Conlon, Fox, Gunby, He, Mubayi, Suk, Verstra\"ete, and Yu~\cite{CFGHHMSVY} conjectured that containment in an iterated blowup of an edge is also necessary, and this was recently confirmed by Gishboliner, Hajebi, and Li~\cite{GHL2026}. While this gives a complete characterization of the polynomial regime for $3$-graphs, determining the precise growth rate within this regime remains a central challenging problem, even for small forbidden hypergraphs.

A natural family of examples is given by the \emph{Berge triangles}, which consist of three distinct edges containing, respectively, the three pairs of a set of three vertices. Up to isomorphism, there are four such $3$-graphs: the tight path $TP_3=\{abc,bcd,cde\}$, the loose triangle $C_3=\{abc,cde,efa\}$, and the $3$-graphs $F_5$ and $K_4^{(3)-}$ defined above. As observed by Mattheus, Mubayi, Nie, and Verstra\"ete~\cite{MMNV}, their off-diagonal Ramsey numbers exhibit remarkably different growth rates. Phelps and R\"odl~\cite{PhelpsRodl} showed that $R(TP_3,K_t^{(3)})$ lies between $\Omega(t^2/\log t)$ and $O(t^2)$. For $F_5$, the upper bound $O(t^3/\log t)$ of Kostochka, Mubayi, and Verstra\"ete~\cite{KMV2014} and the lower bound $\Omega(t^3/(\log t)^4)$ of Mattheus, Mubayi, Nie, and Verstra\"ete~\cite{MMNV} determine the order of growth up to a polylogarithmic factor. In contrast, $R(K_4^{(3)-},K_t^{(3)})=t^{\Theta(t)}$~\cite{ErdosHajnal72,FoxHe}. For $C_3$, Theorem~\ref{thm:Ramsey} gives the first improvement on the upper bound $O(t^{3/2})$ of Kostochka, Mubayi, and Verstra\"ete~\cite{KMV}.

\subsection{Independent sets in hypergraphs}\label{sec:independent-sets}

Theorem~\ref{thm:Ramsey} follows from a lower bound on the independence number of loose-triangle-free hypergraphs (Theorem~\ref{thm:main} below), which we now place in context.
The starting point is the Tur\'an-type bound of Spencer~\cite{Spencer}: every $n$-vertex $r$-graph with average degree $d \ge 1$ satisfies
\begin{equation}\label{eq: spencer}
    \alpha(H)\ge\Big(1-\frac1r\Big)\frac{n}{d^{1/(r-1)}},
\end{equation}
which is tight up to the constant factor for disjoint unions of cliques.
A recurring theme is that local sparsity improves~\eqref{eq: spencer} by a polylogarithmic factor.
For graphs, Ajtai, Koml\'os, and Szemer\'edi~\cite{AKS} showed that triangle-free graphs satisfy $\alpha(G)=\Omega(n\log d/d)$, which Shearer~\cite{Shearer} sharpened to $(1-o(1))\,n\log d/d$.
A standard way to quantify local sparsity is through cycles.
A \emph{Berge cycle} of length $k\ge2$ consists of distinct vertices $v_1,\ldots,v_k$ and distinct edges $e_1,\ldots,e_k$ such that $\{v_i,v_{i+1}\}\subseteq e_i$ for each $i$, where $v_{k+1}=v_1$; the \emph{girth} of a hypergraph is the length of its shortest cycle, and a hypergraph is \emph{linear} if it has no cycle of length two, that is, if any two edges share at most one vertex.
The phenomenon was first exploited for hypergraphs by Koml\'os, Pintz, and Szemer\'edi~\cite{KPS} in their work on Heilbronn's triangle problem, and was established for all uniformities by Ajtai, Koml\'os, Pintz, Spencer, and Szemer\'edi~\cite{AKPSS}: every $n$-vertex \emph{uncrowded} $r$-graph, that is, one of girth at least five, with maximum degree $\Delta$ satisfies
\begin{equation}\label{eq: AKPSS}
    \alpha(H)\ge c_r\, n\left(\frac{\log\Delta}{\Delta}\right)^{1/(r-1)}
\end{equation}
for some constant $c_r>0$; Duke, Lefmann, and R\"odl~\cite{DLR} showed that linearity alone suffices.

Determining the optimal $c_r$ has attracted considerable recent attention~\cite{Iliopoulos,VW}, often via coloring: a \emph{proper coloring} assigns colors to the vertices so that no edge is monochromatic, and since each color class is independent, $\alpha(H)\ge n/\chi(H)$, where the \emph{chromatic number} $\chi(H)$ is the minimum number of colors in a proper coloring.
Dhawan, Methuku, and Vo~\cite{DMV,DMVcoloring} and, independently, Yu and Zhang~\cite{YuZhang,YuZhangList} showed that uncrowded $r$-graphs are list colorable with $(1+o(1))\big((r-1)\Delta/\log\Delta\big)^{1/(r-1)}$ colors, matching the coupon-collector heuristic and confirming a conjecture of Molloy~\cite{MolloyCanaDAM}.
Dhawan, Methuku, and Vo~\cite{DMV} also resolved a conjecture of Verstra\"ete and Wilson~\cite{VW}, showing that for linear hypergraphs the constant $c_r$ can be chosen so that $c_r\to1$ as $r\to\infty$.

A related sparsity condition is triangle-freeness, in the sense of Cooper and Mubayi~\cite{CM} and Li and Postle~\cite{LP}.
A \emph{triangle} consists of three distinct edges $e,f,g$ and three distinct vertices $x,y,z$ such that
\begin{equation}\label{eq: triangle def CM LP}
    \{x,y\}\subseteq e,\qquad \{y,z\}\subseteq f,\qquad
    \{z,x\}\subseteq g,\qquad\text{and}\qquad
    \{x,y,z\}\cap e\cap f\cap g=\varnothing,
\end{equation}
or equivalently, such that each of $(e\cap f)\setminus g$, $(f\cap g)\setminus e$, and $(g\cap e)\setminus f$ is nonempty.
Throughout the paper, \emph{triangle} and \emph{triangle-free} are meant in this sense.
Every loose triangle is a triangle, and the two notions coincide in linear hypergraphs, and in particular in graphs.
For graphs, Johansson~\cite{Johansson} proved that triangle-free graphs of maximum degree $\Delta$ satisfy $\chi(G)=O(\Delta/\log\Delta)$.\footnote{A published account of Johansson's argument is given by Molloy and Reed~\cite{MolloyReed}, who obtain $\chi(G) \leq 160\,\Delta/\log \Delta$ for triangle-free graphs $G$ of sufficiently large maximum degree $\Delta$.}
Frieze and Mubayi~\cite{frieze2013coloring} proved the analogue $\chi(H) = O_r\big((\Delta/\log\Delta)^{1/(r-1)}\big)$ for linear $r$-graphs, and Cooper and Mubayi~\cite[Theorem~2]{CM} replaced linearity by triangle-freeness for $3$-graphs: every triangle-free $3$-graph $H$ with sufficiently large maximum degree $\Delta$ satisfies
\begin{equation}\label{eq: CM coloring}
    \chi(H)\le C\left(\frac{\Delta}{\log\Delta}\right)^{1/2}
\end{equation}
for an absolute constant $C>0$.
Li and Postle~\cite[Theorem~1.4]{LP} extended this by showing that every triangle-free $r$-graph $H$ with sufficiently large maximum degree $\Delta$ satisfies $\chi(H)\le C_r(\Delta/\log\Delta)^{1/(r-1)}$, and hence
$\alpha(H)=\Omega\big(n(\log\Delta/\Delta)^{1/(r-1)}\big)$,
which improves the elementary bound $\Omega(n/\Delta^{1/(r-1)})$ implied by~\eqref{eq: spencer} by a factor of order $(\log\Delta)^{1/(r-1)}$.

The class of hypergraphs without loose triangles is considerably broader: such hypergraphs need not be linear and may contain many triangles, so the preceding bounds do not apply.
For loose-triangle-free $3$-graphs, the best previous bound in terms of $n$ alone is $\alpha(H)=\Omega(n^{2/3})$, which follows from~\eqref{eq: spencer} and the bound $\alpha(H)=\Omega(d)$ of Kostochka, Mubayi, and Verstra\"ete~\cite{KMV}, where $d$ is the average degree.
The stronger bound $\alpha(H)=\Omega(n^{2/3}(\log n)^{1/3})$ is known when $H$ is linear by~\cite[Theorem~1.2]{KMV} and when $H$ is triangle-free by~\eqref{eq: CM coloring}~\cite[Section~1.2]{CM}.
We show that the absence of loose triangles alone already guarantees an improvement over $n^{2/3}$.

\begin{tcolorbox}[enhanced, breakable]
    \begin{theorem}\label{thm:main}
        For every integer $r\ge3$, there is a constant $c_r>0$ such that every $r$-uniform hypergraph $H$ on $n\ge3$ vertices with no loose triangle satisfies
        \[
            \alpha(H)\ge c_r\, n^{2/3}(\log\log n)^{1/3}.
        \]
    \end{theorem}
\end{tcolorbox}

Theorem~\ref{thm:main} implies Theorem~\ref{thm:Ramsey}: for large $t$, let $n\coloneqq\lceil C_r t^{3/2}/\sqrt{\log\log t}\,\rceil$; since $\log\log n = (1+o(1))\log\log t$, Theorem~\ref{thm:main} gives every $C_3$-free $r$-graph on $n$ vertices an independent set of size at least $(1 - o(1))\,c_rC_r^{2/3}\,t\ge t$ once $C_r$ is large.
Conversely, $N\asymp t^{3/2}/\sqrt{\log\log t}$ is equivalent to $t\asymp N^{2/3}(\log\log N)^{1/3}$, so the two theorems are equivalent up to constant factors.
Compared with the linear and triangle-free cases, Theorem~\ref{thm:main} requires neither assumption, at the cost of replacing $\log n$ by $\log\log n$.

Our proof of Theorem~\ref{thm:main} combines probabilistic methods for sampling independent sets with new ideas for controlling triangles involving pairs of vertices of large codegree.
The sampling method for independent sets is inspired by ideas developed in the study of fractional colorings~\cite{Dhawan,MS,DLMVY,DLLMV,YuZhangFractional,bradavc2026coloring}.
In the final step, we apply the theorem of Li and Postle (see Theorem~\ref{thm:LP}). We outline these ideas in Section~\ref{sec:approach}.

\subsection{Our approach}\label{sec:approach}

We now outline the main ideas of the proof; a more detailed overview appears in Section~\ref{sec:overview}.
Recall that a \emph{triangle} is meant in the sense of~\eqref{eq: triangle def CM LP}.
We first describe the argument for $3$-graphs, so let $H$ be a $3$-graph with no loose triangle.

Our strategy is to find a large set of vertices of $H$ that induces a triangle-free subhypergraph of small maximum degree, and to apply the theorem of Cooper and Mubayi (see~\eqref{eq: CM coloring}) to this subhypergraph.
We call a pair of vertices \emph{heavy} if it lies in many edges of $H$, and we let $G$ be the graph whose edges are the heavy pairs.
This graph is triangle-free, since a triangle in $G$ would give a loose triangle in $H$ (see Fig.~\ref{fig: heavy triangle}).
Our aim is to find a large independent set $I$ of $G$ such that all triangles in the induced hypergraph $H[I]$ can be destroyed by deleting few vertices.

To implement this strategy, we combine two key ideas.
First, we show that every triangle in $H[I]$ contains two edges whose union is a four-vertex \emph{$2$-independent} set in $G$, meaning that no two of its vertices are adjacent or have a common neighbor in $G$.
We call these unions \emph{four-vertex configurations}; to destroy all triangles in $H[I]$, it suffices to delete one vertex from each four-vertex configuration contained in $I$.
Second, we show that, for a suitable $\phi>0$, a variant of Johansson's entropy-based method~\cite{Johansson}, as adapted in~\cite{Dhawan} for constructing random independent sets (see also~\cite{MS,bradavc2026coloring}), produces a random independent set $I$ of $G$ that contains each vertex with probability at least $\phi/2$ and each $2$-independent set $S$ with probability at most $\phi^{|S|}$.

Combining these two ideas, each four-vertex configuration is contained in $I$ with probability at most $\phi^4$.
This small probability, together with our bound on the total number of configurations and the fact that every triangle in $H[I]$ contains such a configuration, ensures that the expected number of vertex deletions needed to destroy all triangles is small compared with the expected size of $I$.
Every edge containing no heavy pair is also $2$-independent, so each such edge is contained in $I$ with probability at most $\phi^3$, and this controls the expected number of edges of $H[I]$.
Hence there is an outcome for which deleting one vertex from each four-vertex configuration contained in $I$ leaves a large triangle-free hypergraph of small average degree.
Deleting a small proportion of additional vertices, namely those of large degree, gives the required bound on the maximum degree, allowing us to apply the theorem of Cooper and Mubayi (see~\eqref{eq: CM coloring}).

For general $r$, we combine this argument with a reduction of Kostochka, Mubayi, and Verstra\"ete~\cite{KMV} and some additional ideas, and in the final step we apply the theorem of Li and Postle~\cite{LP} in place of the theorem of Cooper and Mubayi.

\subsection*{Notation}

All logarithms are natural.
For positive functions $f,g$, we write $f=O(g)$ if $f/g$ is bounded above, $f=\Omega(g)$ if $g=O(f)$, $f\asymp g$ or $f=\Theta(g)$ if both hold, $f=o(g)$ if $f/g\to0$, and $f\sim g$ if $f/g\to1$.
A subscript indicates the parameters on which the implicit constants may depend; such parameters are fixed in every limiting statement.
In heuristic discussions only, we write $f\approx g$ for equality up to constant factors, and $f\ll g$ for $f=o(g)$.
For $a \in [0,1]$, an \emph{$a$-random subset} of a set includes each element independently with probability $a$.

For a positive integer $k$, we write $[k]=\{1,\ldots,k\}$.
A \emph{$k$-set} is a set of $k$ vertices, and pairs and triples are $2$-sets and $3$-sets, respectively.
We sometimes abbreviate $\{a,b\}$ to $ab$ and $\{a,b,c\}$ to $abc$.
For a graph $G$, we write $N_G(v)$ for the set of neighbors of $v$, $\deg_G(v)=|N_G(v)|$, and $\Delta(G)$ for the maximum degree of $G$.

Let $H$ be a hypergraph.
We write $e(H)=|E(H)|$.
The \emph{degree} of a vertex is the number of edges containing it.
For $S\subseteq V(H)$, the \emph{induced} hypergraph $H[S]$ has vertex set $S$ and edge set $\{e\in E(H) : e\subseteq S\}$; thus deleting vertices removes the edges that meet the deleted set, rather than shortening them.
We write $E_j(H)$ for the set of edges of size $j$ and $m_j(H)=|E_j(H)|$.
For a set $X\subseteq V(H)$, let
\[
    \dd_j(X,H)=|\{e\in E_j(H):X\subseteq e\}|.
\]
The \emph{$j$-degree} of a vertex $v$ is $\dd_j(v,H) \coloneqq \dd_j(\{v\},H)$, and for $|X|\ge2$, we call $\dd_j(X,H)$ the \emph{$j$-codegree} of $X$.
We write
\[
    \Delta_j(H)=\max_{v\in V(H)}\dd_j(v,H)
    \qquad\text{and}\qquad
    d_j(H)=\frac{j\, m_j(H)}{|V(H)|}
\]
for the maximum and average $j$-degrees, respectively; both are zero if $H$ has no edges of size $j$.
We omit $H$ when it is clear from the context.
A hypergraph is an \emph{antichain} if no edge properly contains another edge.

We write $\pr$ for probability and $\E$ for expectation.
For an event $F$, we write $F^c$ for its complement and $\ind_F$ for its indicator.

\subsection*{Organization}

Section~\ref{sec:overview} gives an overview of the proof, focusing on the $3$-uniform case before describing the additional steps needed for larger uniformities.
In Section~\ref{sec:completion}, we prove Theorem~\ref{thm:main}, following the same sequence of reductions, sampling, and deletions.
Section~\ref{sec:sampling} contains the full definition and analysis of the procedure for sampling random independent sets, which is sketched in Section~\ref{sec:overview-selection}.
Appendix~\ref{app:compression} gives a self-contained proof of the structural result in Lemma~\ref{lem:compression}, following ideas of Kostochka, Mubayi, and Verstra\"ete~\cite{KMV}.

\section{Proof overview}\label{sec:overview}
In this section, we give a detailed overview of our proof.
In Section~\ref{subsec: overview}, we explain the argument in the special case $r = 3$, modulo a technical result, Proposition~\ref{prop: triangle free}.
In Section~\ref{subsec: sampling}, we sketch the proof of Proposition~\ref{prop: triangle free}.
Finally, in Section~\ref{subsec: rank r}, we describe the modifications needed for general $r$.
The full proof appears in Sections~\ref{sec:completion} and~\ref{sec:sampling}.

\subsection{The special case \texorpdfstring{$r = 3$}{r = 3}}\label{subsec: overview}

Let $H$ be an $n$-vertex $3$-graph with no loose triangle, and let $d \coloneqq n^{2/3}(\log\log n)^{1/3}$; we explain this choice below.
Kostochka, Mubayi, and Verstra\"ete~\cite{KMV} showed that $\alpha(H)$ is at least a constant multiple of the average degree of $H$ (see Lemma~\ref{lem:density}), so we may assume that $H$ has average degree at most $d$.
Deleting the vertices of degree greater than $2d$ leaves at least $n/2$ vertices.
For convenience, we ignore this constant-factor loss throughout the overview and assume that $H$ has $n$ vertices and satisfies $\Delta(H)\le 2d$.

Our strategy is to find a large induced subhypergraph of $H$ that is triangle-free in the sense of Cooper and Mubayi~\cite{CM}.
Recall that a \emph{triangle} in a $3$-graph consists of three distinct edges $e, f, g$ and three distinct vertices $a, b, c$ such that $\set{a, b}\subseteq e$, $\set{b, c}\subseteq f$, $\set{a, c}\subseteq g$, and $\set{a, b, c}\cap e\cap f\cap g=\emptyset$.
The Cooper--Mubayi bound~\eqref{eq: CM coloring} implies that every
$m$-vertex triangle-free $3$-graph with maximum degree at most
$\Delta\ge2$ has an independent set of size at least
\begin{equation}\label{eq: triangle free alpha}
    c_0\,m\sqrt{\frac{\log\Delta}{\Delta}}
\end{equation}
for an absolute constant $c_0>0$.
It therefore suffices to find a large set $I \subseteq V(H)$ such that $H[I]$ is triangle-free and has small maximum degree, and then apply~\eqref{eq: triangle free alpha} to obtain the desired independent set in $H$.

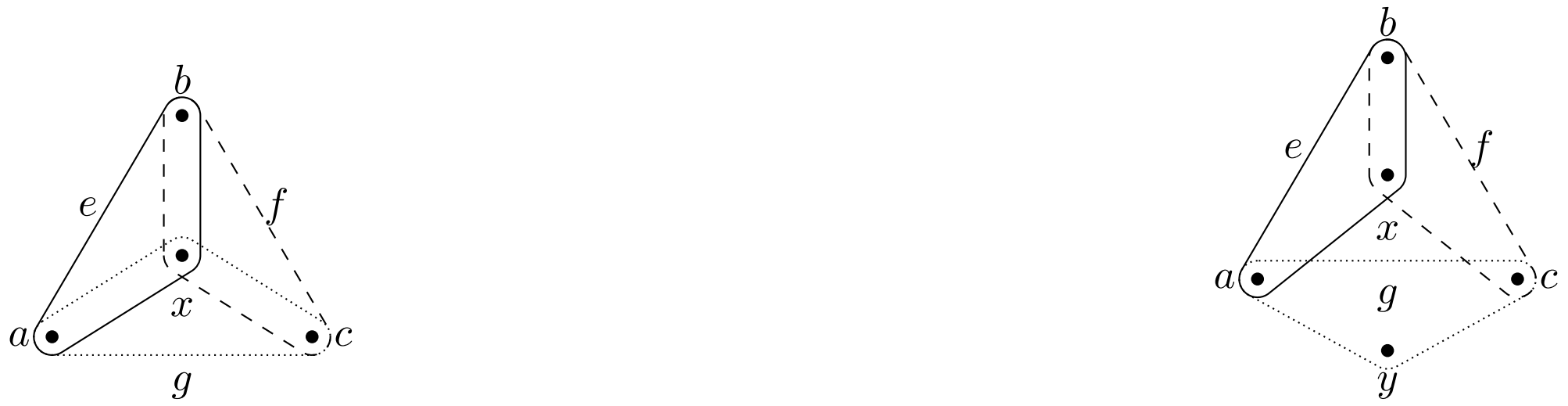
\begin{figure}[htb!]
  \centering
  \begin{subfigure}[b]{0.4\textwidth}
    \centering
    \begin{tikzpicture}[scale=0.55]
      \coordinate (a) at (0,0);
      \coordinate (c) at (4,0);
      \coordinate (b) at (2,3.4);
      \coordinate (x) at (2,1.25);
      \hedge[edgeA]{a}{x}{b}
      \hedge[edgeB]{x}{c}{b}
      \hedge[edgeC]{a}{c}{x}
      \node[vtx, label={left:$a$}] at (a) {};
      \node[vtx, label={above:$b$}] at (b) {};
      \node[vtx, label={right:$c$}] at (c) {};
      \node[vtx, label={[label distance=5pt]270:$x$}] at (x) {};
      \node at (0.55,2.0) {$e$};
      \node at (3.45,2.0) {$f$};
      \node at (2,-0.75) {$g$};
    \end{tikzpicture}
    \caption{$e \cap f \cap g = \set{x}$.}
  \end{subfigure}
  \hfill
  \begin{subfigure}[b]{0.4\textwidth}
    \centering
    \begin{tikzpicture}[scale=0.55]
      \coordinate (a) at (0,0);
      \coordinate (c) at (4,0);
      \coordinate (b) at (2,3.4);
      \coordinate (x) at (2,1.6);
      \coordinate (y) at (2,-1.1);
      \hedge[edgeA]{a}{x}{b}
      \hedge[edgeB]{x}{c}{b}
      \hedge[edgeC]{a}{y}{c}
      \node[vtx, label={left:$a$}] at (a) {};
      \node[vtx, label={above:$b$}] at (b) {};
      \node[vtx, label={right:$c$}] at (c) {};
      \node[vtx, label={[label distance=6pt]270:$x$}] at (x) {};
      \node[vtx, label={below:$y$}] at (y) {};
      \node at (0.55,2.0) {$e$};
      \node at (3.45,2.0) {$f$};
      \node at (2,-0.3) {$g$};
    \end{tikzpicture}
    \caption{$e \cap f = \set{b, x}$ and $y \ne x$.}
  \end{subfigure}
  \caption{The two forms~\eqref{eq: triangle types} of a triangle in a $3$-graph with no loose triangle. In both cases, the edges $e$ and $f$ share the pair $\set{b, x}$.}
  \label{fig: triangle types}
\end{figure}

Because $H$ has no loose triangle, its triangles have a restricted structure.
Indeed, the three edges of a triangle either share a common vertex, or two of them share two vertices.
In either case, after relabeling, the triangle has the form
\begin{equation}\label{eq: triangle types}
    e = \set{a, b, x}, \qquad f = \set{b, c, x}, \qquad g = \set{a, c, y},
\end{equation}
where $a, b, c, x$ are distinct and possibly $y = x$ (see Fig.~\ref{fig: triangle types}).
In particular, the edges $e$ and $f$ share the pair $\set{b, x}$.
A natural first attempt at constructing $I$ is to take a $p$-random subset of $V(H)$ and delete a vertex from each triangle it contains.
This fails, however, when $H$ contains many pairs of large codegree.
Such a pair can play the role of $\set{b, x}$ in many triangles at once, so a $p$-random set may contain far more triangles than vertices.
On the other hand, deleting a single vertex of such a pair destroys all of these triangles simultaneously, which suggests treating pairs of large codegree separately from the remaining triangles.

To this end, we fix a threshold $L \ge 3$ and call a pair of vertices \emph{heavy} if its codegree in $H$ is greater than $L$, and \emph{light} otherwise.
An edge of $H$ is \emph{light} if all three of its pairs are light, and a triangle is \emph{light} if all three of its edges are light.
The \emph{heavy graph} of $H$ is the graph $G$ on $V(H)$ in which $uv \in E(G)$ if and only if $\set{u, v}$ is heavy.
Fix a vertex $v$. Each of the $\deg_G(v)$ heavy pairs containing $v$ lies in more than $L$ edges of $H$, while each edge containing $v$ contains at most two of these pairs.
Since $v$ lies in at most $\Delta(H) \le 2d$ edges of $H$, we obtain $L\deg_G(v)\le 4d$, and hence $\Delta(G)\le\Delta\coloneqq 4d/L$.

For any light triangle in $H$ written as in~\eqref{eq: triangle types}, we call the four-vertex set $e\cup f=\set{a,b,c,x}$ a \emph{bad configuration} (cf.\ Section~\ref{sec:configurations}). If $I$ is independent in $G$, then every triangle of $H[I]$ is light, so deleting one vertex from each bad configuration contained in $I$ destroys every triangle of $H[I]$.
To count bad configurations in $H$, choose a light edge $e$, one of its three pairs, and a distinct light edge $f$ containing that pair.
Since the pair is light, there are at most $L$ choices for $f$.
Every bad configuration arises as $e\cup f$ in this way, so there are at most $3L\,e(H)\le 2ndL=O(ndL)$ bad configurations.
Counting these configurations is more useful than counting light triangles, since a single configuration may arise from many triangles.
We therefore seek a large independent set $I$ of $G$ containing few bad configurations.

A natural attempt is the alteration method: take a $p$-random set $I_0\subseteq V(H)$, then delete one vertex from each edge of $G[I_0]$ and each bad configuration contained in $I_0$.
Since $\E[|I_0|]=np$, improving on the \emph{elementary bound} $\Omega(n/\sqrt d)$ given by Spencer's estimate~\eqref{eq: spencer} requires $p\gg1/\sqrt d$, and hence $dp^2\to\infty$.
The expected number of bad configurations contained in $I_0$ is $O(ndLp^4)$, whereas $\E[|I_0|]=np$.
Hence, to guarantee that deleting these bad configurations costs only a small fraction of $np$, we require $dLp^3$ to be bounded by a sufficiently small constant.
Since $\Delta=4d/L$, we have
\begin{equation}\label{eq: key identity}
    dLp^3=\frac{4(dp^2)^2}{\Delta p}.
\end{equation}
As $dp^2\to\infty$, keeping $dLp^3$ bounded by a sufficiently small constant forces $\Delta p\to\infty$.
However, when $\Delta p\to\infty$, the expected number of edges of $G[I_0]$ can be as large as $n\Delta p^2/2=np\cdot\Delta p/2\gg np$, so we cannot rely on deleting one vertex from each edge of $G[I_0]$ to obtain an independent set of size $\Omega(np)$.

\begin{figure}[htb!]
  \centering
  \begin{tikzpicture}[scale=0.45]
    \coordinate (u1) at (0,2.6);
    \coordinate (u2) at (-2.25,-1.3);
    \coordinate (u3) at (2.25,-1.3);
    \coordinate (x3) at (-3.9,2.25);
    \coordinate (x1) at (0,-4.5);
    \coordinate (x2) at (3.9,2.25);
    \hedge[edgeA]{u1}{x3}{u2}
    \hedge[edgeB]{u2}{x1}{u3}
    \hedge[edgeC]{u3}{x2}{u1}
    \draw[heavy] (u1) -- (u2) -- (u3) -- cycle;
    \node[vtx, label={[label distance=4pt]above:$u_1$}] at (u1) {};
    \node[vtx, label={[label distance=4pt]200:$u_2$}] at (u2) {};
    \node[vtx, label={[label distance=4pt]-20:$u_3$}] at (u3) {};
    \node[vtx, label={[label distance=4pt]left:$x_3$}] at (x3) {};
    \node[vtx, label={[label distance=4pt]below:$x_1$}] at (x1) {};
    \node[vtx, label={[label distance=4pt]right:$x_2$}] at (x2) {};
  \end{tikzpicture}
  \caption{If $u_1, u_2, u_3$ form a triangle in the heavy graph $G$ (thick lines), then the heavy pairs yield edges $\set{u_1, u_2, x_3}$, $\set{u_2, u_3, x_1}$, and $\set{u_1, u_3, x_2}$ that form a loose triangle in $H$.}
  \label{fig: heavy triangle}
\end{figure}
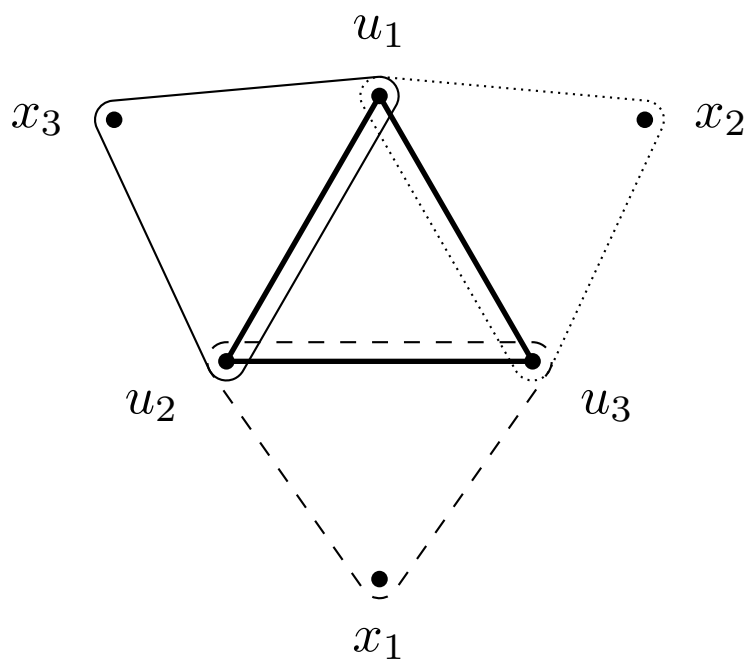

We therefore seek a procedure that samples an independent set $I$ of $G$ directly, while still controlling the bad configurations it contains.
If we could achieve $\E[|I|]=\Omega(np)$ for the values of $p$ considered above, the expected size of $I$ would exceed the greedy bound $n/(\Delta+1)$ by a factor tending to infinity, since $\Delta p\to\infty$.
Fortunately, we can improve on the greedy bound because the heavy graph $G$ is triangle-free.
Indeed, to see why $G$ is triangle-free, suppose that $u_1,u_2,u_3$ form a triangle in $G$.
Since each pair $u_iu_j$ has codegree greater than $L\ge3$, we can successively choose distinct vertices $x_1,x_2,x_3\notin\set{u_1,u_2,u_3}$ such that $\set{u_i,u_j,x_k}\in E(H)$ whenever $\set{i,j,k}=\set{1,2,3}$.
These three edges form a loose triangle in $H$, a contradiction (see Fig.~\ref{fig: heavy triangle}).

Since $G$ is triangle-free, Shearer's theorem~\cite{Shearer} guarantees an independent set of size $\Omega(n\phi)$, where $\phi\approx\log\Delta/\Delta$ and hence $\Delta\phi\approx\log\Delta$.
We seek such an independent set $I$ with the following properties:
\begin{enumerate}[label=(\alph*)]
    \item\label{item: large I} $|I| = \Omega(n\phi)$,
    \item\label{item: few bad configs} $I$ contains $o(n\phi)$ bad configurations, and
    \item\label{item: small degree} all but a sufficiently small constant fraction of the vertices of $I$ have degree $O(d\phi^2)$ in $H[I]$.
\end{enumerate}
Deleting the exceptional vertices in~\ref{item: small degree}, together with one vertex from each bad configuration in~\ref{item: few bad configs}, leaves a set of size $\Omega(n\phi)$ that induces a triangle-free hypergraph of maximum degree $O(d\phi^2)$.

To see that such a set would suffice, suppose for the moment that $H[I]$ resembles the subhypergraph of $H$ induced by a $\phi$-random subset of $V(H)$.
Then $H[I]$ has average degree roughly $d\phi^2$, which gives~\ref{item: small degree} after discarding the few vertices of much larger degree.
Furthermore, as in the alteration argument above, we expect $I$ to contain $O(ndL\phi^4)$ bad configurations, which is $o(n\phi)$ provided that $dL\phi^3 = o(1)$.
Applying~\eqref{eq: triangle free alpha} to the resulting triangle-free hypergraph gives
\[
    \alpha(H) \gtrsim n\phi\sqrt{\frac{\log(d\phi^2)}{d\phi^2}} = n\sqrt{\frac{\log(d\phi^2)}{d}}.
\]
The estimate above improves on the elementary bound $\Omega(n/\sqrt d)$ by a factor of order $\sqrt{\log(d\phi^2)}$.
A convenient choice is $L=\sqrt d/(\log n)^{7/8}$.
Since $\Delta=4d/L$ and $\phi\approx\log\Delta/\Delta$, this gives $\log\Delta\asymp\log n$ and $d\phi^2\asymp(\log n)^{1/4}$.
Consequently, $\log(d\phi^2)\asymp\log\log n$, and the bound above becomes
\[
    \alpha(H)\gtrsim n\sqrt{\frac{\log\log n}{d}},
\]
giving a factor of order $\sqrt{\log\log n}$ over the elementary bound.
This choice also makes the expected number of bad configurations $o(n\phi)$: indeed, by~\eqref{eq: key identity} with $p=\phi$ and $\Delta\phi\approx\log\Delta$,
\[
    dL\phi^3
    \approx\frac{(d\phi^2)^2}{\log\Delta}
    \asymp(\log n)^{-1/2}
    =o(1),
\]
so $O(ndL\phi^4)=o(n\phi)$.
Finally, balancing $n\sqrt{\log\log n/d}$ against the bound $\Omega(d)$ from Lemma~\ref{lem:density} gives $d=n^{2/3}(\log\log n)^{1/3}$ and hence
\[
    \alpha(H)=\Omega\big(n^{2/3}(\log\log n)^{1/3}\big).
\]
Thus the final gain over the bound $\Omega(n^{2/3})$ is $(\log\log n)^{1/3}$.

The difficulty is that a large independent set of $G$ must be chosen according to the structure of $G$, so there is no a priori reason for it to resemble a $\phi$-random set from the point of view of $H$.
One might try to get around this by first taking a $\phi$-random subset $I_0 \subseteq V(H)$, which does behave randomly with respect to $H$, and then passing to a large independent set of $G[I_0]$ using Shearer's bound instead of alteration.
This approach fails as well.
The expected size of $I_0$ is of order $n\phi$, and the expected degree in $G[I_0]$ of a vertex of $I_0$ is at most $\Delta\phi\approx\log\Delta$.
Applying Shearer's bound to $G[I_0]$ therefore only guarantees an independent set of size roughly
\[
    n\phi \cdot \frac{\log(\Delta\phi)}{\Delta\phi} \approx \frac{n\log\log\Delta}{\Delta},
\]
so this estimate falls short of the bound required in~\ref{item: large I} by a factor of order $\log\Delta/\log\log\Delta$.
This loss overwhelms the $\sqrt{\log\log n}$ gain we want to get, and the resulting estimate is in fact worse than the elementary bound $\Omega(n/\sqrt d)$.

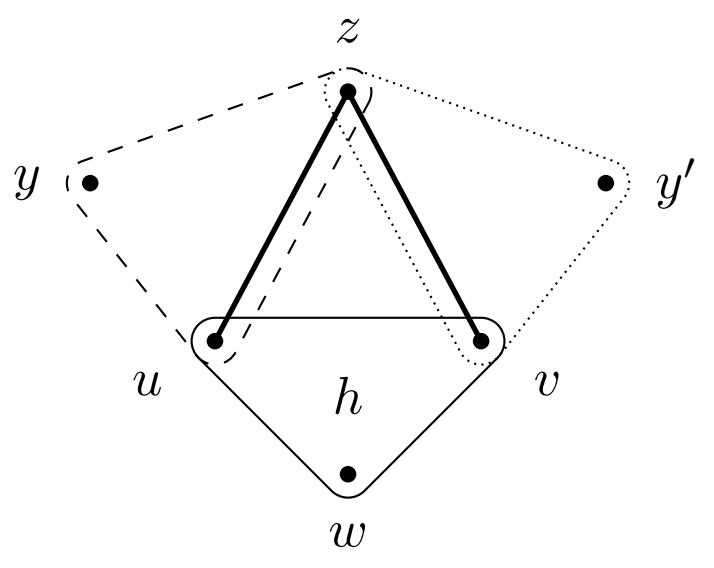
\begin{figure}[htb!]
  \centering
  \begin{tikzpicture}[scale=0.55]
    \coordinate (u) at (0,0);
    \coordinate (v) at (3.2,0);
    \coordinate (w) at (1.6,-1.6);
    \coordinate (z) at (1.6,3.0);
    \coordinate (y) at (-1.5,1.9);
    \coordinate (yp) at (4.7,1.9);
    \hedge[edgeA]{u}{w}{v}
    \hedge[edgeB]{z}{y}{u}
    \hedge[edgeC]{z}{v}{yp}
    \draw[heavy] (z) -- (u);
    \draw[heavy] (z) -- (v);
    \node[vtx, label={[label distance=5pt]200:$u$}] at (u) {};
    \node[vtx, label={[label distance=5pt]-20:$v$}] at (v) {};
    \node[vtx, label={[label distance=4pt]below:$w$}] at (w) {};
    \node[vtx, label={[label distance=4pt]above:$z$}] at (z) {};
    \node[vtx, label={[label distance=4pt]left:$y$}] at (y) {};
    \node[vtx, label={[label distance=4pt]right:$y'$}] at (yp) {};
    \node at (1.6,-0.65) {$h$};
  \end{tikzpicture}
  \caption{If a vertex $z$ is adjacent in $G$ to two vertices $u$ and $v$ of a light edge $h$ (thick lines), then the heavy pairs $zu$ and $zv$ yield edges $\set{z, u, y}$ and $\set{z, v, y'}$ that, together with $h$, form a loose triangle.}
  \label{fig: light edge 2-independent}
\end{figure}

We overcome this obstacle by using the notion of $2$-independence.
A set $S \subseteq V(G)$ is \emph{$2$-independent} in $G$ if any two vertices of $S$ are at a distance greater than $2$ in $G$; equivalently, $S$ is independent and no vertex of $G$ has two neighbors in $S$.
Our key observation is that the sets we need to control, namely light edges and bad configurations, are $2$-independent in $G$.
For a light edge $h$, this holds because a vertex $z$ with two neighbors $u, v \in h$ in $G$ would yield a loose triangle formed by $h$ and two edges through the heavy pairs $zu$ and $zv$ (see Fig.~\ref{fig: light edge 2-independent}).
A bad configuration inherits $2$-independence from the light edges $e$, $f$, and $g$, since each of its six pairs lies in one of them.
It therefore suffices to sample an independent set of $G$ that contains each vertex with probability close to $\phi$, and each $2$-independent set $S$ with probability at most $\phi^{|S|}$, as if its vertices were chosen independently.
This is provided by the following proposition, which yields properties~\ref{item: large I}--\ref{item: small degree} and completes the proof sketch.

\begin{proposition}\label{prop: triangle free}
    Let $G$ be a triangle-free graph with maximum degree at most $\Delta$, and let $0 < \phi < \hp \le 1/2$.
    Then there is a probability distribution over independent sets $I$ of $G$ such that
    \begin{enumerate}[label=(S\arabic*)]
        \item\label{item: vertex prob} $\pr[v \in I] \ge \phi\left(1 - \frac{(1+\hp)\Delta\phi}{\log(\hp/\phi)}\right)$ for every $v \in V(G)$,
        \item\label{item: 2 ISET prob} $\pr[S\subseteq I] \le \phi^{|S|}$ for every $2$-independent set $S \subseteq V(G)$, and
        \item\label{item: all sets prob} $\pr[S\subseteq I]\le\hp^{|S|}$ for every $S\subseteq V(G)$.
    \end{enumerate}
\end{proposition}

For example, if $\phi = \log\Delta/(10\Delta)$ and $\hp = \sqrt{\phi}$, then the error term in~\ref{item: vertex prob} is less than $1/2$ for $\Delta$ sufficiently large, so the expected size of $I$ is within a constant factor of Shearer's bound.
Property~\ref{item: all sets prob} is weaker than~\ref{item: 2 ISET prob} but applies to all sets; it is needed only for $r \ge 4$.

\subsection{Sampling independent sets}\label{sec:overview-selection}\label{subsec: sampling}

We now sketch the proof of Proposition~\ref{prop: triangle free} using Johansson's entropy-based method~\cite{Johansson}, as adapted in~\cite[Section~2]{Dhawan} for constructing random independent sets (see also~\cite{MS} for related work); the details appear in Section~\ref{sec:sampling}.
We process the vertices of $G$ in an arbitrary order $(v_1,\ldots,v_m)$, and write $N_L(v)$ for the set of neighbors of $v$ that precede $v$.
Each vertex $v$ is assigned a weight, initially $p_0(v)=\phi$, and we write $p_i(v)$ for its weight after the $i$-th step.
In the $i$-th step, the vertex $v_i$ is added to $I$ with probability equal to its current weight.
If $v_i$ is added, the weights of its later neighbors are set to $0$, which keeps $I$ independent; otherwise, they are multiplied by $1/(1 - p_{i-1}(v_i))$.
This update rule ensures that, for each vertex $v$, the weights $(p_i(v))_i$ form a martingale.
To keep the weights bounded, a vertex whose weight reaches $\hp$ becomes \emph{bad}: its weight is never updated again, and it is never added to $I$.
By the martingale property, the expected weight of $v_k$ just before it is processed is $\phi$, so
\[
    \pr[v_k \in I] = \phi - \E\big[p_{k-1}(v_k)\,\mathbf{1}\set{v_k \text{ is bad}}\big],
\]
and it remains to show that vertices rarely become bad.

To do so, we track the \emph{entropy contribution} $h_i(v)\coloneqq-p_i(v)\log p_i(v)$.
Since the weights are martingales, a large increase in the weight of $v_k$ is reflected in a decrease of its expected entropy contribution; for $v_k$ to become bad, its weight must grow from $\phi$ to at least $\hp$, which costs entropy of order $\log(\hp/\phi)$ per unit of weight.
It therefore suffices to show that the expected entropy contribution at $v_k$ does not drop much before $v_k$ is processed.
Processing a vertex $v_j \in N_L(v_k)$ decreases it by at most $(1+\hp)\,p_{j-1}(v_j)p_{j-1}(v_k)$ in expectation.
Here the triangle-freeness of $G$ is crucial: since $v_j$ and $v_k$ have no common neighbor, no step changes both of their weights, and so the product $p_i(v_j)p_i(v_k)$ is also a martingale and has expectation $\phi^2$.
Summing over the at most $\Delta$ vertices of $N_L(v_k)$, the expected entropy loss at $v_k$ is at most $(1+\hp)\Delta\phi^2$, and~\ref{item: vertex prob} follows.

The same mechanism gives~\ref{item: 2 ISET prob}.
Fix a $2$-independent set $S$, and for $0 \le i \le m$, consider
\[
    R_i \coloneqq \mathbf{1}_{A_i} \cdot \prod_{\substack{v_k \in S \\ k > i}} p_i(v_k),
\]
where $A_i$ is the event that every vertex of $S$ among $v_1, \ldots, v_i$ was added to $I$.
Intuitively, $R_i$ is the probability that $S \subseteq I$ if the vertices of $S$ that have not yet been processed were included independently with their current weights.
We have $R_0 = \phi^{|S|}$ and $R_m=\mathbf{1}\set{S\subseteq I}$, so it suffices to show that $(R_i)_i$ is a supermartingale.
If $v_i \in S$, then, since $S$ is independent, processing $v_i$ does not change the weights of the other vertices of $S$, and the factor $p_{i-1}(v_i)$ in $R_{i-1}$ is replaced by the indicator that $v_i$ was added, whose conditional expectation is at most $p_{i-1}(v_i)$.
If $v_i \notin S$, then, since $S$ is $2$-independent, $v_i$ has at most one neighbor in $S$, so at most one factor of the product changes, and its conditional expectation is preserved by the martingale property.
Hence $\pr[S \subseteq I] = \E[R_m] \le R_0 = \phi^{|S|}$.
Finally, property~\ref{item: all sets prob} is simpler: a vertex is added only while its weight is below $\hp$, so the probability that all vertices of $S$ are added is at most $\hp^{|S|}$.

\subsection{The extension to larger uniformities}\label{subsec: rank r}

For general $r$, we follow the same strategy after a preliminary reduction of Kostochka, Mubayi, and Verstra\"ete~\cite[Section~3]{KMV} (see Lemma~\ref{lem:compression}), which replaces some edges by smaller subsets while preserving the absence of loose triangles, in such a way that every independent set of the new hypergraph is independent in the original one.
The resulting hypergraph may have edges of several sizes between $3$ and $r$, which is why we prove Theorem~\ref{thm:main2} in this generality, but its larger edges are sparse: it has $O_r(n^2)$ edges of each size at least four, and every triple lies in $O_r(1)$ edges of size four.

The main new difficulty is that two edges may share more than two vertices.
We define heavy pairs using separate thresholds for edges of sizes three and four; as before, the heavy graph $G$ is triangle-free and every light edge is $2$-independent in $G$ (see Lemma~\ref{lem: auxiliaryK3,2ind,covered}).
We then sample $I$ using Proposition~\ref{prop: triangle free} and delete a vertex from each bad configuration contained in $I$ (see Definition~\ref{def: bad configurations} for a formal definition).
The analogues of the configurations from the $3$-uniform case are $2$-independent and are handled as before using~\ref{item: 2 ISET prob}, but two further steps are needed.
First, only light edges of size at least five can survive in $I$, and these are $2$-independent; since the reduction leaves only $O_r(n^2)$ edges of these sizes,~\ref{item: 2 ISET prob} shows that few survive in $I$, and deleting a vertex from each leaves only edges of sizes three and four.
Second, pairs of edges sharing three vertices, as well as the triangles not destroyed by the previous deletions, need not be $2$-independent, so~\ref{item: 2 ISET prob} does not apply to them.
For these, we use the weaker bound~\ref{item: all sets prob}, which applies to all sets but with a larger parameter $\hp$ in place of $\phi$.
The codegree bounds from the reduction of~\cite[Section~3]{KMV} show that these configurations are rare enough (see Lemma~\ref{lem:classification} and Section~\ref{sec:configuration-counts}) for this weaker bound to suffice, at the cost of only a constant factor in $\phi$.
Finally, after deleting vertices of large degree, we apply the coloring theorem of Li and Postle (Theorem~\ref{thm:LP}) to the resulting triangle-free hypergraph with edges of sizes three and four, and with our choice of parameters we again obtain an independent set of size $\Omega_r\big(n^{2/3}(\log\log n)^{1/3}\big)$.

\section{Proof of Theorem~\ref{thm:main}}\label{sec:completion}

In this section, we prove the following extension of Theorem~\ref{thm:main} to hypergraphs in which every edge has size between $3$ and $r$.
\begin{theorem}
\label{thm:main2}
    For every integer $r\ge3$, there is a constant $c_r>0$ such that every hypergraph $\cH$ on $n\ge3$ vertices with no loose triangle and with all edge sizes between $3$ and $r$ satisfies
    \[
        \alpha(\cH)\ge c_r\, n^{2/3}(\log\log n)^{1/3}.
    \]
\end{theorem}
The proof follows the same outline as the $3$-uniform case in Section~\ref{subsec: overview}, with two additional complications: edges may have several sizes, and pairs of edges may now intersect in three or more vertices.
In Section~\ref{sec:reduction}, we recall the results of Kostochka, Mubayi, and Verstra\"ete~\cite{KMV} that we need.
In Section~\ref{sec:sparse-case}, we use these results to handle the dense case and to reduce the maximum degrees, and we state the coloring theorem of Li and Postle~\cite{LP}.
In Section~\ref{sec:auxiliary}, we define the heavy graph and show that it is triangle-free and that light edges are $2$-independent in it.
In Sections~\ref{sec:configurations} and~\ref{sec:configuration-counts}, we identify the configurations whose deletion destroys every triangle, and we bound the number of configurations of each type.
Finally, in Section~\ref{sec:completion-5}, we choose the parameters, bound the cost of the deletions, and complete the proof.

\subsection{Preliminaries}\label{sec:reduction}
We first recall a few structural results from~\cite[Section~3]{KMV}.
Throughout this subsection, all hypergraphs have edge sizes between $3$ and a fixed integer $r\ge3$.
The term ``triangle'' in~\cite{KMV} refers to a loose triangle in our terminology.
\begin{definition}[Expandable sets]\label{def: expandable}
    A nonempty set $X\subseteq V(\cH)$ is \emph{expandable in $\cH$} if, for every set $Y\subseteq V(\cH)\setminus X$ with $|Y|\le2r$, there is an edge of $\cH$ containing $X$ and disjoint from $Y$. 
\end{definition}
Every edge of $\cH$ is expandable in $\cH$, and every expandable set is contained in an edge and hence has size at most $r$.
As shown in~\cite[Lemma~3.2]{KMV}, loose triangles formed by expandable sets can be ``expanded'' to loose triangles formed by edges. We prove it here for completeness.
\begin{lemma}\label{lem:expand}
    If $\cH$ contains no loose triangle, then no three expandable sets in $\cH$ form a loose triangle.
\end{lemma}

\begin{proof}
    Suppose that expandable sets $X_1,X_2,X_3$ form a loose triangle.
    Choose an edge $e_1\supseteq X_1$ disjoint from $(X_2\cup X_3)\setminus X_1$, then an edge $e_2\supseteq X_2$ disjoint from $(e_1\cup X_3)\setminus X_2$, and finally an edge $e_3\supseteq X_3$ disjoint from $(e_1\cup e_2)\setminus X_3$. The edges $e_1,e_2,e_3$ form a loose triangle in $\cH$, a contradiction.
\end{proof}

Recall that $m_j(\cH)$ denotes the number of edges of size $j$ in a hypergraph $\cH$, that $\dd_j(X,\cH)$ denotes the number of these edges containing $X$, and that $\cH$ is an antichain if no edge of $\cH$ properly contains another.
The following lemma collects the facts we need from the proofs of~\cite[Lemmas~3.2--3.4, Corollary~3.5, and Lemma~3.6]{KMV}. We give a proof in Appendix~\ref{app:compression} for completeness.

\begin{lemma}\label{lem:compression}
    For every integer $r\ge3$, there is a constant $C_r>0$ such that the following holds. If $\cH$ is an $n$-vertex hypergraph with no loose triangle in which every edge has size between $3$ and $r$, then there is a hypergraph $\cH'$ on $V(\cH)$ with the following properties:
    \begin{enumerate}[label=\textup{(\roman*)}]
        \item\label{item: compression structure} $\cH'$ has no loose triangle, has edge sizes between $3$ and $r$, and is an antichain;
        \item\label{item: compression independent} every independent set of $\cH'$ is independent in $\cH$;
        \item\label{item: compression codegree} every triple $X\subseteq V(\cH')$ satisfies $\dd_4(X,\cH')\le 2r$; and
        \item\label{item: compression sparse} $m_j(\cH')\le C_r n^2$ for every $4\le j\le r$.
    \end{enumerate}
\end{lemma}

We also use the following bound, which is obtained in~\cite[Section~3.2, Cases~1 and~2]{KMV}.

\begin{lemma}\label{lem:density}
    If $\cH$ has no loose triangle and $d_3=3m_3(\cH)/|V(\cH)|$, then
    \[
        \alpha(\cH)\ge\frac{2d_3}{9(r+1)}.
    \]
\end{lemma}
\subsection{Reduction to the sparse case}\label{sec:sparse-case}\label{sec:degree-reduction}
Let $\cH$ be an $n$-vertex hypergraph as in Theorem~\ref{thm:main2}, and set
\begin{equation}\label{eq: def of d}
    d\coloneqq n^{2/3}(\log \log n)^{1/3}.
\end{equation}
All statements below hold for $n$ sufficiently large in terms of $r$.
Let $\cH'$ be the hypergraph obtained by applying Lemma~\ref{lem:compression} to $\cH$.
If $d_3(\cH') \ge d$, then Lemma~\ref{lem:density} together with Lemma~\ref{lem:compression}\ref{item: compression independent} yields
\[
    \alpha(\cH)\ge\alpha(\cH')\ge\frac{2d}{9(r+1)},
\]
which proves Theorem~\ref{thm:main2}, as desired.
We may therefore assume that $d_3(\cH')<d$.
Since every independent set of $\cH'$ is independent in $\cH$, it suffices to find a large independent set in $\cH'$.

We first pass to an induced subhypergraph $\cG$, which has bounded maximum degrees for edges of each uniformity. By Lemma~\ref{lem:compression}\ref{item: compression sparse}, we have  $d_j(\cH') = j\,m_j(\cH')/n \le rC_r n$ for every $4 \le j \le r$. 
For each $j\in\{3,\ldots,r\}$ with $d_j(\cH')>0$, we delete every vertex $v$ with $\dd_j(v,\cH')>2(r-2)d_j(\cH')$. Since the $j$-degrees sum to $n\,d_j(\cH')$, at most $n/(2(r-2))$ vertices are deleted for each $j$, and so at most $n/2$ vertices are deleted in total.
Let $\cG$ be the subhypergraph of $\cH'$ induced by the remaining vertices, and let $n_0\coloneqq|V(\cG)|$, so that $n/2\le n_0\le n$.

Since $\cG$ is an induced subhypergraph of $\cH'$, it is an antichain with no loose triangle by Lemma~\ref{lem:compression}\ref{item: compression structure}. Setting $b_r\coloneqq \max\{2r, 2r(r-2)C_r\}$, we have
\begin{align}
    \Delta_3(\cG)&\le b_r d,\label{eq:degree3}\\
    \Delta_j(\cG)&\le b_r n\qquad\text{for } 4\le j\le r,\label{eq:degreehigher}\\
    \dd_4(X,\cG)&\le b_r\qquad\text{for every triple }X\subseteq V(\cG),\label{eq:degree34}
\end{align}
where~\eqref{eq:degree34} follows from Lemma~\ref{lem:compression}\ref{item: compression codegree}.
In particular, \eqref{eq:degree3} and \eqref{eq:degreehigher} imply that $m_3(\cG) \le b_r nd/3$ and $m_j(\cG) \le b_r n^2/j$ for every $4 \le j \le r$. The above steps are summarized in the following chart:
\[
    \cH\xrightarrow{\text{Lemma~\ref{lem:compression}}}\cH'
    \xrightarrow{\text{high degree vertex deletion}}\cG.
\]
The rest of the proof focuses on showing that $\alpha(\cG)\ge c_rd$, since every independent set of $\cG$ is independent in $\cH'$, and hence in $\cH$. Before we proceed further, we state the following result of Li and Postle on coloring triangle-free hypergraphs with edges of sizes $3$ and $4$.
\begin{theorem}[{Li--Postle~\cite[Theorem 1.4]{LP}}]\label{thm:LP}
    There is an absolute constant $C$ with the following property.
    Let $\cK$ be a triangle-free hypergraph in which every edge has size $3$ or $4$.
    Let $D_3,D_4\ge8$ be such that, for each $j\in\{3,4\}$, every vertex of $\cK$ lies in at most $D_j$ edges of size $j$.
    Then
    \[
        \chi(\cK)\le C\max\left\{
        \left(\frac{D_3}{\log D_3}\right)^{1/2},
        \left(\frac{D_4}{\log D_4}\right)^{1/3}
        \right\}.
    \]
\end{theorem}

\subsection{Heavy graphs and \texorpdfstring{$2$}{2}-independent sets}\label{sec:auxiliary}

In Section~\ref{sec:deletions}, we will show that deleting one vertex from each edge of size at least five contained in our random independent set costs only a small fraction of the sampled vertices. We now thus define heavy pairs using only edges of sizes three and four, and use two thresholds $L_3$ and $L_4$ as functions of $d$ (see \eqref{eq: def of d}) and $n$. Let   
\begin{equation}
\label{eq: L3L4}
     L_3\coloneqq\frac{\sqrt d}{(\log n)^{7/8}} = \frac{n^{1/3}(\log\log n)^{1/6}}{(\log n)^{7/8}},\qquad
    L_4\coloneqq \frac{n L_3}{d} = \frac{n^{2/3}}{(\log\log n)^{1/6}\,(\log n)^{7/8}}. 
\end{equation}
These choices will be useful in the final computations in Section~\ref{sec:completion-5}. For now, it is sufficient to observe that since $n\to \infty$ and $r$ is fixed,  $L_3$ and $L_4$ satisfy
\begin{equation}\label{eq:threshold-relations}
    L_3>2r,\qquad L_4>2rb_r,\qquad L_4=\frac nd L_3,
    \qquad\text{and}\qquad L_j\le n^{(j-2)/3}\quad\text{for } j\in\set{3,4}.
\end{equation}
\begin{definition}[Heavy Graph]
    \label{dfn: heavygraph}
    We say that a pair of vertices $P\subseteq V(\cG)$ is \emph{heavy} if $\dd_3(P,\cG)>L_3$ or $\dd_4(P,\cG)>L_4$, and \emph{light} otherwise.
An edge of $\cG$ is \emph{light} if all of its underlying pairs are light.
Let the \emph{heavy graph} $G$ denote the graph with vertex set $V(\cG)$ whose edges are the heavy pairs. 
\end{definition}
The use of heavy and light pairs, with different codegree thresholds, also appears in the work of Kostochka, Mubayi, and Verstra\"ete~\cite[Section~3.2]{KMV} and M\'eroueh~\cite[Section~3]{Meroueh}.

We let $\cG_0\subseteq \cG$ be the hypergraph formed by the light edges of $\cG$. The reductions so far can be summarized as follows:
\[
    \cH\xrightarrow{\text{Lemma~\ref{lem:compression}}}\cH'
    \xrightarrow{\text{vertex deletion}}\cG
    \xrightarrow{\text{only light edges}}\cG_0.
\]
Moreover, observe that if $I$ is an independent set in $G$, then every edge of $\cG[I]$ is light, and so $\cG[I]=\cG_0[I]$. First we show that $G$ has bounded degree. 
\begin{lemma}
    \label{lem: aux degree}
    We have $\Delta(G) \le \Delta$, where $\Delta \coloneqq 5b_rd/L_3$. 
\end{lemma}
\begin{proof}
    Fix $v\in V(\cG)$. Note that, by counting the edges of size 3 in $\cG$ adjacent to $v$, we have: 
    $$ L_3 \cdot  |\{u: \dd_3(\{u,v\})> L_3\}|\le \sum_{u\in N_G(v)}\dd_3(\{u,v\},\cG) \le 2\dd_3(v,\cG)\le 2\Delta_3(\cG).$$
    Similarly, for edges of $\cG$ of size 4 adjacent to $v$, we have
     $$ L_4 \cdot  |\{u: \dd_4(\{u,v\})> L_4\}|\le \sum_{u\in N_G(v)}\dd_4(\{u,v\},\cG) \le 3\Delta_4(\cG).$$
     Note that by definition of $G$, any $u$ adjacent to $v$ either has $\dd_3(\{u,v\},\cG)> L_3$ or $\dd_4(\{u,v\},\cG)> L_4$. Consequently, the two inequalities above together with the bounds in \eqref{eq:degree3}, \eqref{eq:degreehigher}, and~\eqref{eq:threshold-relations} give
    \[
       \deg_G(v)\le \frac{2\Delta_3(\cG)}{L_3}+\frac{3\Delta_4(\cG)}{L_4}
        \le\frac{2b_rd}{L_3}+\frac{3b_rn}{L_4}=\frac{5b_rd}{L_3}.\qedhere
    \]
\end{proof}
We now prove crucial structural properties of $G$. 
\begin{lemma}
\label{lem: expandableheavypairs}
    Every pair $\{x,y\}\in E(G)$ is expandable in $\cG$. 
\end{lemma}
\begin{proof}
   Let $P$ be a heavy pair. Fix $Y\subseteq V(\cG)\setminus P$ with $|Y|\le2r$. 
   \begin{itemize}
       \item  If $\dd_3(P,\cG)>L_3$, then the number of edges of size 3 that contain $P$ and intersect $Y$ is at most $|Y| \le 2r < L_3$. So there is an edge containing $P$ that is disjoint from $Y$. 
       \item   If $\dd_4(P,\cG)>L_4$, then in view of~\eqref{eq:degree34}, $\dd_4(P\cup \{y\}, \cG)\le b_r$ for each vertex $y\in Y$. So $Y$ intersects with at most $2rb_r<L_4$ edges of size 4 that contain $P$. 
   \end{itemize}
    In either case, some edge containing $P$ is disjoint from $Y$, so $P$ is expandable. 
\end{proof}
As a corollary of the above lemma, we obtain the following.
\begin{lemma}
    \label{lem: auxiliaryK3,2ind,covered}
   The graph $G$ is triangle-free. Moreover, 
    \begin{enumerate}[label=(\roman *)]
        \item\label{item: aux 2ind} every edge of $\cG_0$ is $2$-independent in $G$; and
        \item\label{item: aux covered} any set $U \subseteq V(\cG)$ is $2$-independent in $G$ if every pair of vertices of $U$ lies in an edge of $\cG_0$.
    \end{enumerate}
\end{lemma}
\begin{proof}
    Since each edge in $G$ is expandable in $\cG$, $G$ is triangle-free in view of Lemma~\ref{lem:expand}. 
    
    Let $e\in E(\cG_0)$. Since every pair of vertices in $e$ is light, $e$ is independent in $G$. Assume for the sake of contradiction that two vertices $x,y\in e$ have a common neighbor $z$ in $G$ with $z\notin e$.
    The sets $\{x,z\}$, $\{y,z\}$ are expandable in $\cG$ by Lemma~\ref{lem: expandableheavypairs}. Further, $e$ is expandable in $\cG$ since it is an edge. These 3 expandable sets of $\cG$ form a loose triangle contradicting Lemma~\ref{lem:expand}. Thus $e$ is $2$-independent.

The vertices of $U$ are pairwise nonadjacent in $G$, since every pair of $U$ lies in a light edge.
    If some vertex $v$ had two neighbors $p, q \in U$ in $G$, then $v$ would have two neighbors in an edge of $\cG_0$ containing $\set{p, q}$, contradicting~\ref{item: aux 2ind}.
\end{proof}

\subsection{The bad configurations}\label{sec:configurations}
Let $I$ be an independent set of $G$; in Section~\ref{sec:completion-5}, $I$ will be the random independent set given by Proposition~\ref{prop: triangle free}.
Then every edge of $\cG[I]$ is light, and our goal is to delete a small number of vertices from $I$ so that the remaining vertices induce a triangle-free hypergraph with edges of sizes three and four only, to which we can apply Theorem~\ref{thm:LP}.

In the $3$-uniform case, every light triangle contains two edges that share a pair of vertices and whose union is $2$-independent in $G$, so it suffices to delete a vertex from each such union contained in $I$.
For general $r$, we first destroy every edge of size at least five contained in $I$ by deleting a vertex from each. Among the remaining triangles, some can be destroyed by deleting a vertex from the union of two edges, either because they share two vertices and their union is $2$-independent in $G$, or because they share three vertices. We count the remaining triangles separately.

Let $\tilde{\cG}$ be the hypergraph consisting of the edges of $\cG_0$ of sizes three and four.
The reductions from the hypergraph $\cH$ in Theorem~\ref{thm:main2} can be summarized as follows:
\[
    \cH\xrightarrow{\text{Lemma~\ref{lem:compression}}}\cH'
    \xrightarrow{\text{vertex deletion}}\cG
    \xrightarrow{\text{light edges}}\cG_0
    \xrightarrow{\text{sizes }3,4}\tilde{\cG}.
\]
\begin{definition}[Bad configurations]\label{def: bad configurations}
    Let $\cP$ be the set of unordered pairs $\{e,f\}$ of distinct edges of $\tilde{\cG}$ such that $|e\cap f|=2$ and $e\cup f$ is $2$-independent in $G$. We define the following four types of bad configurations.
    \begin{enumerate}[label=\textup{(B\arabic*)}]
        \item\label{item: conf large} \emph{Large edges.} Let $\cE_{\ge5}$ be the set of edges of $\cG_0$ of size at least five.
        \item\label{item: conf B} \emph{Unions of pairs in $\cP$.} For $s\in\{4,5,6\}$, let
        \[
            \cB_s\coloneqq\{e\cup f:\{e,f\}\in\cP,\ |e\cup f|=s\}
            \qquad\text{and}\qquad \cB\coloneqq\cB_4\cup\cB_5\cup\cB_6.
        \]
        \item\label{item: conf C} \emph{Edges sharing a triple.} Let $\cC = \{e\cup f: |e|=|f|=4,|e\cap f|=3, e,f\in \tilde{\cG}\}$, i.e., the set of vertex sets of pairs of edges of size 4 sharing 3 vertices. 
        \item\label{item: conf T}\label{dfn: lighttrianglesk} \emph{Uncovered triangles.}
        A triangle $T=\{e,f,g\}$ in $\tilde{\cG}$ is \emph{uncovered} if every two of its edges intersect in at most two vertices and 
        $$\{\{e,f\},\{f,g\},\{e,g\}\}\cap \cP = \emptyset.$$
        Let $\cT$ be the family of all uncovered triangles, and write $V(T)\coloneqq e\cup f\cup g$ for $T=\{e,f,g\}\in\cT$.
    \end{enumerate}
\end{definition}

We now show that deleting one vertex from each of these configurations leaves a triangle-free hypergraph whose edges have sizes three and four.
For a family $\cA$ of vertex sets, we write $\cA(I)\coloneqq\{X\in\cA:X\subseteq I\}$, and write $\cT(I)\coloneqq\{T\in\cT:V(T)\subseteq I\}$. The following observation justifies the name ``bad configurations'' in the above definition. 
\begin{obs}\label{obs: deletions}
    Let $I$ be an independent set of $G$.
    Delete one vertex from each member of $\cE_{\ge5}(I)$, $\cB(I)$, and $\cC(I)$, and one vertex from each triangle in $\cT(I)$.
    The remaining vertices induce in $\cG$ a triangle-free hypergraph whose edges have sizes three and four.
\end{obs}

\begin{proof}
    Since $I$ is independent in $G$, every edge of $\cG[I]$ is light.
    Every such edge of size at least five belongs to $\cE_{\ge5}(I)$
    and is removed when its chosen vertex is deleted.
    Thus every remaining edge has size three or four and belongs
    to $\tilde{\cG}[I]$.

    Now consider any triangle $T$ in $\tilde{\cG}[I]$.
    If two of its edges share at least three vertices, then, since
    $\tilde{\cG}$ is an antichain with edge sizes three and four,
    both have size four and share exactly three vertices.
    Their union therefore belongs to $\cC(I)$.

    Otherwise, every two edges of $T$ intersect in at most two vertices.
    If some pair of edges of $T$ belongs to $\cP$, its union belongs
    to $\cB(I)$.
    If no pair belongs to $\cP$, then $T$ is uncovered by definition,
    so $T\in\cT(I)$.

    In each case, a vertex of $V(T)$ is deleted, so at least one edge
of $T$ is removed.
Thus the remaining vertices of $I$ induce in $\cG$ a triangle-free
hypergraph whose edges have sizes three and four, as desired.
\end{proof}

The following lemma describes the structure of uncovered triangles which we use later in Section~\ref{sec:configuration-counts}.
\begin{lemma}\label{lem:classification}
    Every uncovered triangle admits a labeling $T=\{e,f,g\}$ such that
    \[
        |e|=4\qquad\text{and}\qquad |e\setminus(f\cup g)|=1.
    \]
\end{lemma}

\begin{proof}
    Fix an uncovered triangle $T=\{e,f,g\}$ in $\tilde{\cG}$.
    We will show that, up to relabeling, these edges satisfy the properties above.
    If all three pairwise intersections had size one, then $T$ would be a loose triangle, contradicting the assumption that $\tilde{\cG}$ contains no loose triangle.
    Therefore, at least one pair of edges, say $\{e,f\}$, satisfies $|e\cap f|=2$.

    We make the following two observations.
    \begin{enumerate}[label=(\roman*), ref=(\roman*)]
        \item\label{it: obsa}
        \textit{Whenever two edges $e,f\in T$ share two vertices,
        there is a pair of vertices in $e\cup f$ that lies in none
        of the three edges of $T$.}
        Indeed, otherwise
        Lemma~\ref{lem: auxiliaryK3,2ind,covered}\ref{item: aux covered}
        would imply that $e\cup f$ is $2$-independent in $G$,
        giving $\{e,f\}\in\cP$, contrary to the assumption that
        $T$ is uncovered.

        \item\label{it: obsb}
        \textit{Two edges of size three in $T$ cannot share two vertices.}
        To see this, let us assume the edges have vertex sets
        $\{a,b,c\}$ and $\{a,b,d\}$.
        The third edge of $T$ contains both $c$ and $d$, since it meets
        each of these two edges outside the other.
        Every pair in $\{a,b,c,d\}$ then lies in one of the three
        edges, again giving $\{\{a,b,c\},\{a,b,d\}\}\in\cP$,
        a contradiction (see Fig.~\ref{fig: triangle types}).
    \end{enumerate}

    By~\ref{it: obsb} and the fact that $|e\cap f|=2$, the edges $e,f$ cannot both have size three,
    so, after interchanging their labels if necessary, we may assume that $|e|=4$.
    Note that $e\cap g$ is not contained in $e\cap f$, for otherwise
    $(e\cap g)\setminus f$ would be empty, contrary to the definition of a triangle.
    Hence
    \[
        3\le |e\cap(f\cup g)|
        = |(e\cap f)\cup(e\cap g)|\le4.
    \]
    If this intersection has size three, then $e$ has exactly one vertex outside $f\cup g$, as required.

    Suppose now that the above intersection has size four.
    In this case, we have an uncovered triangle $T=\{e,f,g\}$ with $|e|=4$ and
    \[
        |e\cap g|=2,\qquad |e\cap f|=2,\qquad
        (e\cap f)\cap(e\cap g)=\emptyset,\qquad
        e\subseteq f\cup g.
    \]
    Note that the third equality above also implies that $e\cap(f\cap g)=\emptyset.$ 
    
    We now consider two cases: $|f|=3$ or $|g|=3$
    (these are symmetric, so we consider only $|f|=3$),
    and $|f|=|g|=4$.
    \begin{itemize}
        \item If $|f|=3$, we must have $|f\cap g|=1$ since
        $e\cap(f\cap g)=\emptyset$.
        The situation then corresponds to Fig.~\ref{fig: classification a}
        (here $g$ may have four vertices, but the fourth vertex is omitted).
        However, this contradicts~\ref{it: obsa}, since every pair of vertices
        in $e\cup f$ is contained in one of $e,f,g$.

\begin{figure}[b!]
  \centering
  \begin{subfigure}[b]{0.31\textwidth}
    \centering
    \begin{tikzpicture}[scale=0.55]
      \coordinate (q) at (0,1.6);
      \coordinate (p) at (1.2,2.0);
      \coordinate (t) at (0,-0.6);
      \coordinate (s) at (1.2,-1.0);
      \coordinate (c) at (4.2,0.5);
      \hedgefour[edgeA]{q}{t}{s}{p}
      \hedge[edgeB]{q}{c}{p}
      \hedge[edgeC]{t}{s}{c}
      \foreach \v in {q,t,p,s,c} \node[vtx] at (\v) {};
      \node at (-0.6,0.5) {$e$};
      \node at (2.9, 2) {$f$};
      \node at (2.9,-0.85) {$g$};
    \end{tikzpicture}
    \caption{$|f|=3$.}
    \label{fig: classification a}
  \end{subfigure}
  \hfill
  \begin{subfigure}[b]{0.31\textwidth}
    \centering
    \begin{tikzpicture}[scale=0.55]
      \coordinate (a1) at (-0.6,2.2);
      \coordinate (a2) at (0.6,2.2);
      \coordinate (b1) at (-2.0,-0.3);
      \coordinate (b2) at (-1.3,-1.4);
      \coordinate (c1) at (2.0,-0.3);
      \coordinate (c2) at (1.3,-1.4);
      \hedgefour[edgeA]{a1}{b1}{b2}{a2}
      \hedgefour[edgeB]{a1}{c2}{c1}{a2}
      \hedgefour[edgeC]{b1}{b2}{c2}{c1}
      \foreach \v in {a1,a2,b1,b2,c1,c2} \node[vtx] at (\v) {};
      \node at (-1.75,1.2) {$e$};
      \node at (1.75,1.5) {$f$};
      \node at (0,-2.15) {$g$};
    \end{tikzpicture}
    \caption{$|f|=|g|=4$, $|f\cap g|=2$.}
    \label{fig: classification b}
  \end{subfigure}
  \hfill
  \begin{subfigure}[b]{0.31\textwidth}
    \centering
    \begin{tikzpicture}[scale=0.55]
      \coordinate (q) at (0,1.6);
      \coordinate (p) at (1.2,2.0);
      \coordinate (t) at (0,-0.6);
      \coordinate (s) at (1.2,-1.0);
      \coordinate (c) at (4.2,0.5);
      \coordinate (u) at (3.4,2.5);
      \coordinate (v) at (3.4,-1.5);
      \hedgefour[edgeA]{q}{t}{s}{p}
      \hedgefour[edgeB]{q}{c}{u}{p}
      \hedgefour[edgeC]{t}{s}{v}{c}
      \foreach \w in {q,t,p,s,c,u,v} \node[vtx] at (\w) {};
      \node at (-0.65,0.5) {$e$};
      \node at (2.6, 3) {$f$};
      \node at (2.6,-1.95) {$g$};
    \end{tikzpicture}
    \caption{$|f|=|g|=4$, $|f\cap g|=1$.}
    \label{fig: classification c}
  \end{subfigure}
  \caption{In~\subref{fig: classification a} and~\subref{fig: classification b},
every pair in $e\cup f$ lies in one of $e,f,g$, contradicting~\ref{it: obsa}.
In~\subref{fig: classification c}, we swap the labels $e$ and $f$ to obtain the desired labeling.}
  \label{fig: classification}
\end{figure}
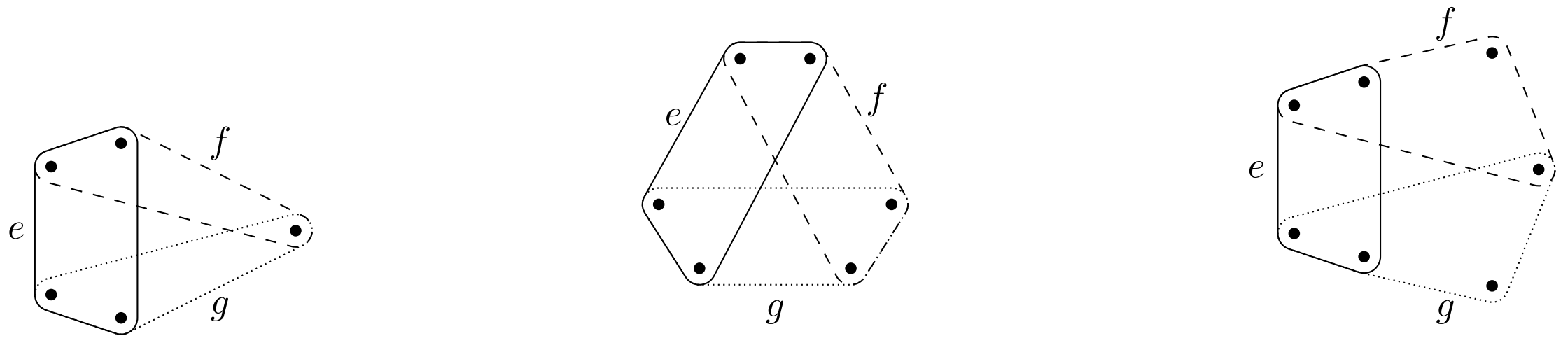

        \item If $|f|=|g|=4$, then we need to consider two subcases:
        \begin{itemize}[label=$*$]
            \item If $|f\cap g|=2$, then the situation corresponds to
            Fig.~\ref{fig: classification b}.
            In this case, once again every pair of vertices in $e\cup f$
            is contained in one of $e,f,g$, contradicting~\ref{it: obsa}.

            \item If $|f\cap g|=1$, then the situation corresponds to
            Fig.~\ref{fig: classification c}.
            In this case, we have
            \[
                |f|=4,\qquad |f\setminus(e\cup g)|=1,
            \]
            so swapping the labels $e$ and $f$ completes the proof. \qedhere
        \end{itemize}
    \end{itemize}
\end{proof}

\subsection{Counting bad configurations}\label{sec:configuration-counts}

In this subsection, we show that there are few bad configurations of each type. In Section~\ref{sec:deletions}, we combine the following counting estimates with the sampling bounds to show that all four families can be destroyed by deleting a small number of vertices.

First we consider~\ref{item: conf large}. For $5 \le j \le r$, the number of members of $\cE_{\ge 5}$ of size $j$ is at most $m_j(\cG)$, so by~\eqref{eq:degreehigher},
\begin{equation}\label{eq:E-count}
    |\set{X \in \cE_{\ge 5} : |X| = j}| \le m_j(\cG) \le \frac{b_r n^2}{j} \qquad \text{for every } 5 \le j \le r.
\end{equation}
When $r \le 4$, the family $\cE_{\ge 5}$ is empty.

Next we consider~\ref{item: conf B}.
The family $\cB_4$ includes the configurations from the
$3$-uniform case, while $\cB_5$ and $\cB_6$ are their analogues
involving edges of size four.
Since every edge of $\tilde{\cG}$ is light, we have
$\dd_j(P,\tilde{\cG})\le L_j$ for every pair $P$ and
$j\in\{3,4\}$.
Moreover, double counting gives
\[
    \sum_P\dd_j(P,\tilde{\cG})
    =\binom j2 m_j(\tilde{\cG})
    \le \frac{j-1}{2}\,n\Delta_j(\cG)
    \qquad\text{for }j\in\{3,4\}.
\]
Counting pairs of edges that share a pair of vertices therefore yields
\begin{equation}\label{eq:B}
\begin{aligned}
    |\cB_4|
    &\le\sum_P\binom{\dd_3(P,\tilde{\cG})}2
    \le\frac{L_3}2\sum_P\dd_3(P,\tilde{\cG})
    \overset{\eqref{eq:degree3}}{\le}\frac{b_r}2\,ndL_3,\\
    |\cB_5|
    &\le\sum_P\dd_3(P,\tilde{\cG})\,\dd_4(P,\tilde{\cG})
    \le L_4\sum_P\dd_3(P,\tilde{\cG})
    \overset{\eqref{eq:degree3}}{\le}b_r\,ndL_4,
    \\
    |\cB_6|
    &\le\sum_P\binom{\dd_4(P,\tilde{\cG})}2
    \le\frac{L_4}2\sum_P\dd_4(P,\tilde{\cG})
    \overset{\eqref{eq:degreehigher}}{\le}\frac{3b_r}4\,n^2L_4.
\end{aligned}
\end{equation}
To bound $|\cC|$, we choose the first edge, one of its four triples, and a second edge containing this triple; by~\eqref{eq:degreehigher} and~\eqref{eq:degree34}, this gives
\begin{equation}\label{eq:C-count}
    |\cC|\le4b_r\, m_4(\tilde{\cG})\le b_r^2n^2.
\end{equation}
Finally, we claim that for every $n^{-1/3}\le q\le1$, 
\begin{equation}\label{eq:triangle-count}
    \sum_{T\in\cT}q^{|V(T)|}
    \le1680\,b_r^3 n^3q^8(\log\log n)^{2/3}.
\end{equation}
By~\eqref{eq: def of d}, \eqref{eq:degree3}, and~\eqref{eq:degreehigher}, for $i\in\{3,4\}$ we have
\begin{equation}\label{eq:edgecounts-reduced}
    m_i(\tilde{\cG})\le\frac{n\Delta_i(\cG)}i
    \le b_rn^{(i+2)/3}(\log\log n)^{1/3}.
\end{equation}
Similarly, for $j\in\{3,4\}$, $s\in\{1,2\}$, and an $s$-set $X$,
\[
    \dd_j(X,\tilde{\cG})\le b_rn^{(j-s)/3}(\log\log n)^{1/3}.
\]
For $s=1$, this follows from~\eqref{eq:degree3} and~\eqref{eq:degreehigher}; for $s=2$, it follows from $\dd_j(X,\tilde{\cG})\le L_j\le n^{(j-2)/3}$.

Label each $T\in\cT$ as in Lemma~\ref{lem:classification}, and put $i\coloneqq|f|$, $j\coloneqq|g|$, and $s\coloneqq|f\cap g|$.
Then $i,j\in\{3,4\}$, since edges of size at least five were deleted in forming $\tilde\cG$, and $s\in\{1,2\}$, since $f\cap g\neq\varnothing$ in a triangle and $|f\cap g|\ge3$ would give $f\cup g\in\cC$, contradicting that $T$ is uncovered.
Since $|e\setminus(f\cup g)|=1$, we have
\[
    |V(T)|=|e\setminus (f\cup g)| + |f\cup g| =1+ (i+j-s).
\]
Fix $(i,j,s)$. We count the triangles with these parameters by choosing $f$, then $g$, then $e$.
\begin{itemize}
    \item \emph{Choosing $f$.} Since $f$ is an edge of size $i$ in $\tilde\cG$, by~\eqref{eq:edgecounts-reduced} there are at most $m_i(\tilde\cG)\le b_rn^{(i+2)/3}(\log\log n)^{1/3}$ choices.
    \item \emph{Choosing $g$.} The intersection $X\coloneqq f\cap g$ is one of at most $\binom{i}{s}\le\binom42=6$ subsets of $f$ of size $s$, and $g$ is an edge of size $j$ containing $X$, so there are at most $\dd_j(X,\tilde\cG)$ choices for $g$ given $X$.
    If $s=1$, then $\dd_j(X,\tilde\cG)\le\Delta_j(\cG)$, which is at most $b_rd=b_rn^{2/3}(\log\log n)^{1/3}$ for $j=3$ by~\eqref{eq:degree3} and at most $b_rn$ for $j=4$ by~\eqref{eq:degreehigher}.
    If $s=2$, then $X$ is a pair contained in an edge of $\tilde\cG$, and hence light, so $\dd_j(X,\tilde\cG)\le L_j\le n^{(j-2)/3}$ by~\eqref{eq:threshold-relations}.
    In all cases, there are at most $b_rn^{(j-s)/3}(\log\log n)^{1/3}$ choices for $g$ given $X$.
    \item \emph{Choosing $e$.} Since $|e|=4$ and $|e\setminus(f\cup g)|=1$, the set $e\cap(f\cup g)$ is one of at most $\binom{|f\cup g|}{3}\le\binom73=35$ triples in $f\cup g$, and each such triple lies in at most $b_r$ edges of size four by~\eqref{eq:degree34}.
\end{itemize}

Writing $k=i+j-s+1\le8$, the contribution of $q^{|V(T)|}$ in the sum in \eqref{eq:triangle-count} from all $T\in \cT$ with a particular choice of $(i,j,s)$ is at most
\[
    210\,b_r^3(\log\log n)^{2/3}n^{(k+1)/3}q^k.
\]
Since $n^{1/3}q\ge1$, we have
\[
    n^{(k+1)/3}q^k
    =n^{1/3}(n^{1/3}q)^k
    \le n^{1/3}(n^{1/3}q)^8
    =n^3q^8.
\]
Summing over the at most eight choices of $(i,j,s)$ proves~\eqref{eq:triangle-count}.

\subsection{Completing the proof}\label{sec:completion-5}\label{sec:parameters}\label{sec:deletions}
We now combine the results of the previous subsections.
We first choose the parameters and sample an independent set $I$ of the heavy graph using Proposition~\ref{prop: triangle free}, which we restate for convenience.

\begin{proposition*}[Restatement of Proposition~\ref{prop: triangle free}]
    Let $G$ be a triangle-free graph with maximum degree at most $\Delta$, and let $0 < \phi < \hp \le 1/2$.
    Then there is a probability distribution over independent sets $I$ of $G$ such that
    \begin{enumerate}[label=\textup{(S\arabic*)}]
        \item $\pr[v \in I] \ge \phi\left(1 - \frac{(1+\hp)\Delta\phi}{\log(\hp/\phi)}\right)$ for every $v \in V(G)$,
        \item $\pr[S\subseteq I] \le \phi^{|S|}$ for every $2$-independent set $S \subseteq V(G)$, and
        \item $\pr[S\subseteq I]\le\hp^{|S|}$ for every $S\subseteq V(G)$.
    \end{enumerate}
\end{proposition*}

We then show that deleting a vertex from each bad configuration in Definition~\ref{def: bad configurations} contained in $I$, together with the vertices of large degree, costs only a small fraction of $I$.
Finally, we apply Theorem~\ref{thm:LP} to the resulting triangle-free hypergraph.

\subsubsection*{Choosing the parameters}
Recall that $d = n^{2/3}(\log\log n)^{1/3}$ from~\eqref{eq: def of d} and $L_3, L_4$ were defined in~\eqref{eq: L3L4} as follows: 
\begin{equation}
\label{eq: L3L42}
     L_3\coloneqq\frac{\sqrt d}{(\log n)^{7/8}} = \frac{n^{1/3}(\log\log n)^{1/6}}{(\log n)^{7/8}},\qquad
    L_4\coloneqq \frac{n L_3}{d} = \frac{n^{2/3}}{(\log\log n)^{1/6}\,(\log n)^{7/8}}. 
\end{equation}
We now set
\begin{equation}\label{eq:parameters}
\begin{gathered}
    \Delta\coloneqq\frac{5b_rd}{L_3} = 5b_rn^{1/3}(\log\log n)^{1/6} (\log n)^{7/8},\\[0.3cm]
    \phi\coloneqq\frac{\log\Delta}{100\Delta},\qquad\text{and}\qquad
    \hp\coloneqq n^{-3/10}.
\end{gathered}
\end{equation}
We need to control the expected sizes of $\cC(I)$ and the family of uncovered triangles $\cT(I)$. The sets in these families need not be $2$-independent and are therefore handled using~\ref{item: all sets prob}. 

We estimate $\phi$ as a function of $n$. We have $\log\Delta=\frac13\log n+O(\log\log n)$, so $\log \Delta \sim \frac13 \log n$ and
\begin{equation}\label{eq:alpha-order}
    \phi=\frac{\log\Delta}{500b_r\sqrt d\,(\log n)^{7/8}}
    =\Theta_r\!\left(n^{-1/3}(\log n)^{1/8}(\log \log n)^{-1/6}\right).
\end{equation}
In particular, $0<\phi<\hp\le1/2$ and $\log(\hp/\phi)\sim\frac1{30}\log n$, so
\[
    \frac{(1+\hp)\Delta\phi}{\log(\hp/\phi)} = \frac{(1+\hp)\log\Delta}{100\log(\hp/\phi)}\longrightarrow\frac1{10}.
\]
By Lemmas~\ref{lem: aux degree} and~\ref{lem: auxiliaryK3,2ind,covered}, we have $\Delta(G)\le \Delta$ and $G$ is triangle-free, so we may apply Proposition~\ref{prop: triangle free} to $G$ with these values of $\Delta$, $\phi$, and $\hp$.
This gives a random independent set $I$ of $G$ satisfying~\ref{item: 2 ISET prob} and~\ref{item: all sets prob} such that $\pr[v\in I]\ge\phi/2$ for every vertex $v\in V(G)$. Recall that $|V(\cG)|=|V(G)|=n_0\ge n/2$, and so it follows that
\begin{equation}\label{eq:sample-size}
    \E[|I|]\ge\frac{n_0\phi}{2}\ge\frac{n\phi}{4}.
\end{equation}

\subsubsection*{The cost of the deletions}
We bound the expected number of bad configurations of each type contained in $I$, in the order of Definition~\ref{def: bad configurations}, beginning with the $2$-independent configurations. 

First we deal with the edges in $\cE_{\ge 5}$ that are contained in the random independent set $I$. Note that in view of Lemma~\ref{lem: auxiliaryK3,2ind,covered}\ref{item: aux 2ind}, every edge of $\cG_0$ is 2-independent. Using the upper bound on $\cE_{\ge 5}$ in \eqref{eq:E-count} and the estimate on $\phi$ in~\eqref{eq:alpha-order} together with the upper bound in~\ref{item: 2 ISET prob}, we obtain that
\begin{equation}\label{eq:large-edge-loss}
    \E[|\cE_{\ge5}(I)|]
    \le\sum_{j=5}^r \frac{b_rn^2}{j}\,\phi^j
    =O_r(n^2\phi^5)=n\phi\cdot O_r(n\phi^4)=o(n\phi).
\end{equation}
Next, we bound the expected number of vertex sets in $\cB$ that are contained in $I$. Recall that $\cB$ contains sets of vertices of the form $e\cup f$ where $\{e,f\}\in \cP$, and thus by definition, $e\cup f$ is 2-independent.  Thus in view of~\ref{item: 2 ISET prob} and the upper bounds on $\cB_s$ in~\eqref{eq:B}, we have
\begin{equation}\label{eq:B-loss-bound}
    \E[|\cB(I)|]
    \le\frac{b_r}2\, ndL_3\phi^4+b_r\, ndL_4\phi^5+\frac{3b_r}4\, n^2L_4\phi^6
    = n\phi \cdot O_r\big(dL_3\phi^3 + dL_4\phi^4 + nL_4\phi^5\big).
\end{equation}
We now upper bound the terms on the right-hand side of the above equation.  
Using the identity $d = n^{2/3} (\log\log n)^{1/3}$, the estimate on $\phi$ in~\eqref{eq:alpha-order}, and the estimates on $L_3$ and $L_4$ in \eqref{eq: L3L42}, we have that
\begin{align}
    dL_3\phi^3&=O_r\big((\log n)^{-1/2}\big),\label{eq:dL-alpha3}\\
    dL_4\phi^4&
    =O_r\big((\log n)^{-3/8}(\log \log n)^{-1/2}\big),\label{eq:dL4-alpha4}\\
    nL_4\phi^5
    &=O_r\big((\log n)^{-1/4}(\log \log n)^{-1}\big).\label{eq:nL4-alpha5}
\end{align}
All three quantities tend to zero, and so we have
\begin{align}
\label{eq:B-loss}
    \E[|\cB(I)|]=o(n\phi).
\end{align}
    We now bound $\E[|\cC(I)|]$. The vertex sets in $\cC$ need not be $2$-independent. However, $|\cC|$ is small enough that the weaker bound in~\ref{item: all sets prob} is sufficient. In view of the upper bound on $|\cC|$ in~\eqref{eq:C-count} and the fact that $n\phi \ge n^{2/3}$,
\begin{equation}\label{eq:C-loss}
    \E[|\cC(I)|]\le b_r^2n^2\hp^5=b_r^2n^{1/2}=o(n\phi).
\end{equation}
Now we bound $\E[|\cT(I)|]$. As with the sets in $\cC$, the unions of uncovered triangles need not be $2$-independent, so we use the bound in~\ref{item: all sets prob} together with the weighted count in~\eqref{eq:triangle-count}.
Since $\hp=n^{-3/10}\ge n^{-1/3}$, applying~\eqref{eq:triangle-count} with $q=\hp$ gives
\begin{equation}\label{eq:T-loss}
    \E[|\cT(I)|]
    \le\sum_{T\in\cT}\hp^{|V(T)|}
    \le1680\,b_r^3(\log\log n)^{2/3}n^3\hp^8
    =O_r\big(n^{3/5}(\log\log n)^{2/3}\big)
    =o(n\phi).
\end{equation}
Here, and in~\eqref{eq:C-loss}, we use $n\phi\ge n^{2/3}$ for sufficiently large $n$, which follows from~\eqref{eq:alpha-order}.

Finally, we control the degrees.
Let
\begin{equation}\label{eq:final-degrees}
    D_3\coloneqq20b_rd\phi^2\qquad\text{and}\qquad D_4\coloneqq20b_rn\phi^3,
\end{equation}
and for $j\in\{3,4\}$, let
\[
    Z_j\coloneqq\{v\in I:\dd_j(v,\cG[I])\ge D_j\}.
\]
Summing $\dd_j(v,\cG[I])$ over $v\in Z_j$, we must have 
\begin{equation}
\label{eq: Zjbound}
    |Z_j|\le jm_j(\cG[I])/D_j.
\end{equation}
Since $I$ is independent in $G$, we have $\cG[I]=\cG_0[I]$.
Every edge of $\cG_0$ is $2$-independent in $G$ by Lemma~\ref{lem: auxiliaryK3,2ind,covered}\ref{item: aux 2ind}.
In view of~\ref{item: 2 ISET prob}, we then have
\[
    \E[m_j(\cG[I])]
    =\sum_{\substack{e\in E(\cG_0)\\ |e|=j}}\pr[e\subseteq I]
    \le m_j(\cG_0)\phi^j
    \le m_j(\cG)\phi^j.
\]
Taking expectations on both sides in \eqref{eq: Zjbound}, and using that  
$j\,m_j(\cG)\le |V(\cG)|\Delta_j(\cG)\le n\Delta_j(\cG)$, we obtain
\[
    \E[|Z_j|] \le \frac{\E[j\,m_j(\cG[I])]}{D_j}
    \le\frac{j\,m_j(\cG)\phi^j}{D_j}
    \le\frac{n\Delta_j(\cG)\phi^j}{D_j}.
\]
Substituting~\eqref{eq:degree3}, \eqref{eq:degreehigher}, and~\eqref{eq:final-degrees} yields
\begin{align*}
    \E[|Z_3|]
    &\le\frac{nb_rd\phi^3}{20b_rd\phi^2}
    =\frac{n\phi}{20},\\
    \E[|Z_4|]
    &\le\frac{nb_rn\phi^4}{20b_rn\phi^3}
    =\frac{n\phi}{20}.
\end{align*}
Consequently,
\begin{equation}\label{eq:high-degree-loss}
    \E[|Z_3\cup Z_4|]
    \le\E[|Z_3|]+\E[|Z_4|]
    \le\frac{n\phi}{10}.
\end{equation}
Finally, we note that by~\eqref{eq:parameters} and $\log\Delta\sim\frac13\log n$,
\begin{equation}\label{eq:d-alpha2}
    d\phi^2=\frac{(\log\Delta)^2}{(500b_r)^2(\log n)^{7/4}}
    =\Theta_r\big((\log n)^{1/4}\big).
\end{equation}

\subsubsection*{Final alteration and putting everything together}
We now delete a vertex from each set in $\cE_{\ge 5}(I)$, $\cB(I)$, and $\cC(I)$, and from the vertex set of each triangle in $\cT(I)$. We further delete every vertex in $Z_3\cup Z_4$. 
Let $I'\subseteq I$ be the set of remaining vertices, and set $\cG'\coloneqq\cG[I']$. 
By Observation~\ref{obs: deletions}, $\cG[I']$ is a triangle-free hypergraph whose edges have sizes three and four.
Finally, since $I'\subseteq I\setminus(Z_3\cup Z_4)$, the maximum $j$-degree of $\cG[I']$ is less than $D_j$ for $j\in\set{3,4}$.

By~\eqref{eq:large-edge-loss}--\eqref{eq:T-loss}, the expected number of vertices deleted for bad configurations is $o(n\phi)$, and in particular at most $n\phi/20$.
Together with~\eqref{eq:sample-size} and~\eqref{eq:high-degree-loss}, this gives
\[
    \E[|I'|]\ge n\phi\left(\frac14-\frac1{10}-\frac1{20}\right)=\frac{n\phi}{10}.
\]
Fix an outcome with $|I'|\ge n\phi/10$.
By~\eqref{eq:d-alpha2} and the identity $n\phi^3 = (d\phi^2)^{3/2}/(\log\log n)^{1/2}$, we have
\[
    D_3=\Theta_r\big((\log n)^{1/4}\big)\qquad\text{and}\qquad
    D_4=\Theta_r\big((\log n)^{3/8}(\log \log n)^{-1/2}\big),
\]
so both tend to infinity, with $\log D_3\sim \frac14\log \log n$ and $\log D_4\sim\frac38\log \log n$.
Theorem~\ref{thm:LP} therefore gives
\begin{align*}
    \chi(\cG')
    &\le C\max\left\{\left(\frac{D_3}{\log D_3}\right)^{1/2},
    \left(\frac{D_4}{\log D_4}\right)^{1/3}\right\}\\
    &\le C'_r\,\phi\max\left\{\left(\frac {d}{\log \log n}\right)^{1/2},\left(\frac {n}{\log \log n}\right)^{1/3}\right\}
    =C'_r\,\phi\, n^{1/3}(\log \log n)^{-1/3}
\end{align*}
for a constant $C'_r>0$ depending only on $r$, where the last equality uses~\eqref{eq: def of d}; note that both terms in the maximum are equal.
Hence some color class of $\cG'$ has size at least
\[
    \frac{|I'|}{\chi(\cG')}\ge\frac1{10C'_r}\,n^{2/3}(\log \log n)^{1/3}.
\]
This color class is independent in $\cG$, hence in $\cH'$, and therefore in $\cH$.
Taking $c_r\coloneqq\min\{2/(9(r+1)),1/(10C'_r)\}$ covers both this case and the dense case from Section~\ref{sec:sparse-case}, proving the desired bound for all sufficiently large $n$.
Since $\alpha(\cH)\ge1$ for $n\ge3$, decreasing $c_r$ if necessary covers the finitely many remaining integers $n\ge3$, proving Theorem~\ref{thm:main2} and hence Theorem~\ref{thm:main}.

\section{The random independent set procedure}\label{sec:sampling}

In this section, we prove Proposition~\ref{prop: triangle free}, which we restate for convenience.

\begin{proposition*}[Restatement of Proposition~\ref{prop: triangle free}]
    Let $G$ be a triangle-free graph with maximum degree at most $\Delta$, and let $0 < \phi < \hp \le 1/2$.
    Then there is a probability distribution over independent sets $I$ of $G$ such that
    \begin{enumerate}[label=\textup{(S\arabic*)}]
        \item $\pr[v \in I] \ge \phi\left(1 - \frac{(1+\hp)\Delta\phi}{\log(\hp/\phi)}\right)$ for every $v \in V(G)$,
        \item $\pr[S\subseteq I] \le \phi^{|S|}$ for every $2$-independent set $S \subseteq V(G)$, and
        \item $\pr[S\subseteq I]\le\hp^{|S|}$ for every $S\subseteq V(G)$.
    \end{enumerate}
\end{proposition*}

The procedure and its analysis use Johansson's entropy-based method~\cite{Johansson}, as adapted in~\cite[Section~2]{Dhawan} for constructing random independent sets (see also the related work of Martinsson and Steiner~\cite{MS}).
Throughout this section, $G$ is a triangle-free graph on $m$ vertices with maximum degree at most $\Delta$, and $0<\phi<\hp\le1/2$.

\subsection{The procedure}

We process the vertices of $G$ in an arbitrary order $(v_1, \ldots, v_m)$.
For each vertex $v$, let $N_L(v)$ and $N_R(v)$ denote the sets of neighbors of $v$ that precede and follow $v$ in this ordering, respectively, and write $d_L(v) \coloneqq |N_L(v)|$.

The procedure maintains a weight $p_i(v)$ for each vertex $v$, where $i$ indexes the steps; initially, every vertex has weight $\phi$.
When a vertex $v_i$ is processed, it is added to $I$ with probability equal to its current weight.
If $v_i$ is added, the weights of its later neighbors are set to $0$; otherwise, they are scaled up so that each weight remains a martingale.
To keep the weights bounded, a vertex whose weight reaches $\hp$ becomes \emph{bad}: its weight is never updated again, and it is never added to $I$.
The set $B_i$ consists of the vertices that are bad after the $i$-th step.
The precise procedure is given in Algorithm~\ref{algorithm: fcp k3}, and we let $I$ denote its output.
\vspace{10pt}
\begin{breakablealgorithm}
\caption{Random Independent Set Procedure for Triangle-Free Graphs}\label{algorithm: fcp k3}
\begin{flushleft}
\textbf{Input}: An $m$-vertex triangle-free graph $G$, an ordering $(v_1, \ldots, v_m)$ of $V(G)$, and parameters $0 < \phi < \hp \le 1/2$. \\
\textbf{Output}: An independent set $I$ in $G$.
\end{flushleft}

\begin{itemize}
\item \textbf{Initialize:} Set $I = B_0 = \emptyset$.
For each $i \in [m]$, set $p_0(v_i) = \phi$.

\item \textbf{Iterate:} For $i = 1, \ldots, m$, do the following:
\begin{enumerate}[itemsep=3pt, label=(\arabic{enumi}), ref=(\arabic{enumi})]
\item\label{step: bad skip} If $v_i\in B_{i-1}$, set $a_i=0$ and $p_i(v_j)=p_{i-1}(v_j)$ for all $j$.
\item If $v_i\notin B_{i-1}$, do the following:
\begin{enumerate}[itemsep=3pt, label=(\alph{enumii}), ref=(\arabic{enumi})(\alph{enumii})]
\item\label{step: no update} For each $j$ such that $v_j \notin N_R(v_i) \setminus B_{i-1}$, set $p_i(v_j) = p_{i-1}(v_j)$.
\item\label{step: coin flip} Sample $a_{i} \sim \Ber(p_{i-1}(v_i))$. If $a_{i} = 1$, add $v_i$ to $I$.
\item For each $v_j \in N_R(v_i)\setminus B_{i-1}$, do the following:
\begin{enumerate}[itemsep=3pt, label=(\roman{enumiii}), ref=(\arabic{enumi})(\alph{enumii})(\roman{enumiii})]
\item\label{step: set zero} if $a_{i} = 1$, set $p_i(v_j) = 0$;
\item\label{step: weight boost} if $a_{i} = 0$, set $p_i(v_j) = \frac{p_{i-1}(v_j)}{1 - p_{i-1}(v_i)}$.
\end{enumerate}
\end{enumerate}
\item\label{step: bad update} Set $B_i=\{v_j\,:\,p_i(v_j)\ge \hp\}$.
\end{enumerate}
\end{itemize}
\end{breakablealgorithm}
\vspace{10pt}

By Steps~\ref{step: set zero} and~\ref{step: bad update}, the output $I$ is an independent set.
Let $\cF_i$ denote the information determined by $a_1,\ldots,a_i$, so that conditioning on $\cF_i$ means conditioning on all choices made in the first $i$ steps; in particular, $\cF_0$ is trivial.
Every weight is at most $\hp/(1-\hp)\le 1$: an update multiplies a weight below $\hp$ by a factor of at most $1/(1-\hp)$, and the weight of a vertex that becomes bad is never updated again.
In particular, the Bernoulli trials in Step~\ref{step: coin flip} are well defined.
Moreover, every weight is either $0$ or at least $\phi$, since each update either sets a weight to $0$, leaves it unchanged, or increases it.

Our first lemma shows that, for each fixed $j$, the weights $p_i(v_j)$ form a martingale in $i$.

\begin{lemma}\label{lemma: K3 martingale}
    For all $1 \le i < j \le m$, we have $\E[p_i(v_j) \mid \cF_{i-1}] = p_{i-1}(v_j)$.
    In particular, $\E[p_i(v_j)] = \phi$ for all $0 \le i < j$.
\end{lemma}

\begin{proof}
    Fix $j \in [m]$ and $i \in [j-1]$.
    If $v_i \in B_{i-1}$, then $p_i(v_j) = p_{i-1}(v_j)$ by Step~\ref{step: bad skip}, so we may assume that $v_i \notin B_{i-1}$.
    If $v_i \notin N_L(v_j)$, then $v_j \notin N_R(v_i)$, and so $p_i(v_j) = p_{i-1}(v_j)$ by Step~\ref{step: no update}.
    Similarly, if $v_j \in B_{i-1}$, then $v_j \notin N_R(v_i) \setminus B_{i-1}$, and again $p_i(v_j) = p_{i-1}(v_j)$ by Step~\ref{step: no update}.
    We may therefore assume that $v_i \in N_L(v_j)$ and $v_j \notin B_{i-1}$.
    Then, by Steps~\ref{step: set zero} and~\ref{step: weight boost},
    \[
        \E[p_i(v_j) \mid \cF_{i-1}]
        = p_{i-1}(v_i) \cdot 0 + (1 - p_{i-1}(v_i)) \cdot \frac{p_{i-1}(v_j)}{1 - p_{i-1}(v_i)} = p_{i-1}(v_j).
    \]
    The second statement follows by taking expectations, since $p_0(v_j) = \phi$.
\end{proof}

\subsection{The probability of containing a given set}

With Lemma~\ref{lemma: K3 martingale} in hand, we prove properties~\ref{item: 2 ISET prob} and~\ref{item: all sets prob}.

\begin{lemma}\label{lemma: 2 ISET prob}
    Let $S$ be a $2$-independent set in $G$.
    Then $\pr[S \subseteq I] \le \phi^{|S|}$.
\end{lemma}

\begin{proof}
    We may assume that $S$ is nonempty, and write $S = \{v_{i_1}, \ldots, v_{i_s}\}$, where $s=|S|$ and $i_1<\cdots<i_s$.
    For $0 \le i \le m$, let $S_{\le i} \coloneqq S \cap \{v_1, \ldots, v_i\}$ and $S_{> i} \coloneqq S \setminus S_{\le i}$, and define $S_{<i}$ and $S_{\ge i}$ analogously.
    Consider the random variables
    \[
        R_i \coloneqq \mathbf{1}\set{S_{\le i}\subseteq I}\prod_{u \in S_{> i}}p_{i}(u), \qquad 0 \le i \le m.
    \]
    Intuitively, $R_i$ is the probability that $S \subseteq I$ if the unprocessed vertices of $S$ were included independently with their current weights.
    Since whether $v_k \in I$ is decided in the $k$-th step, $R_i$ is determined by $\cF_i$.
    Moreover, $R_0 = \phi^s$, since every vertex has initial weight $\phi$, and $R_{i_s} = \mathbf{1}\set{S \subseteq I}$.
    It therefore suffices to show that $(R_i)_i$ is a supermartingale, i.e., that $\E[R_i\mid \cF_{i-1}] \le R_{i-1}$ for each $i \in [m]$.

    Fix $i \in [m]$, and suppose first that $v_i \in S$.
    Since $S$ is independent, $v_i \notin N_L(u)$ for every $u \in S_{> i}$, and so $p_i(u) = p_{i-1}(u)$ for all such $u$.
    Since $v_i$ is added to $I$ with probability $p_{i-1}(v_i)$ if it is not bad, and with probability $0$ otherwise, we have
    \begin{align*}
        \E[R_i\mid \cF_{i-1}] &= \prod_{u \in S_{> i}}p_{i-1}(u)\cdot\E[\mathbf{1}\set{S_{\le i}\subseteq I}\mid \cF_{i-1}] \\
        &= \mathbf{1}\set{S_{< i}\subseteq I}\Big(\prod_{u \in S_{> i}}p_{i-1}(u)\Big) p_{i-1}(v_i)\,\mathbf{1}\set{v_i \notin B_{i-1}} \\
        &\le \mathbf{1}\set{S_{< i}\subseteq I}\prod_{u \in S_{\ge i}}p_{i-1}(u) = R_{i-1}.
    \end{align*}

    Now suppose that $v_i \notin S$, so that $S_{\le i} = S_{<i}$ and $S_{>i} = S_{\ge i}$.
    If $v_i \notin N_L(u)$ for every $u \in S_{> i}$, then $R_i = R_{i-1}$ and there is nothing to prove.
    Otherwise, since $S$ is $2$-independent, $v_i$ has exactly one neighbor $u^*$ in $S$, and $u^*\in S_{>i}$; the weights of all other vertices in $S_{>i}$ are unchanged in the $i$-th step.
    Hence, by Lemma~\ref{lemma: K3 martingale},
    \[
        \E[R_i\mid \cF_{i-1}]
        = \mathbf{1}\set{S_{< i}\subseteq I}\prod_{u \in S_{> i}\setminus\set{u^*}}p_{i-1}(u)\cdot\E[p_{i}(u^*)\mid \cF_{i-1}]
        = R_{i-1}.
    \]

    Taking expectations gives $\E[R_i]\le\E[R_{i-1}]$ for each $i$, and hence $\pr[S\subseteq I]=\E[R_{i_s}]\le R_0=\phi^s$.
\end{proof}

Property~\ref{item: all sets prob} holds for all sets $S$ and has a simpler proof.

\begin{lemma}\label{lemma: all sets prob}
    For every set $S \subseteq V(G)$, we have $\pr[S \subseteq I] \le \hp^{|S|}$.
\end{lemma}

\begin{proof}
    We may assume that $S$ is nonempty, and write $S=\{v_{i_1},\ldots,v_{i_s}\}$ with $i_1<\cdots<i_s$.
    For $0 \le j \le s$, let $E_j$ be the event that $v_{i_1},\ldots,v_{i_j}$ are all added to $I$, so that $E_0$ is the certain event.
    Conditional on $\cF_{i_j-1}$, the vertex $v_{i_j}$ is added to $I$ with probability less than $\hp$: a vertex that is not bad has weight below $\hp$, and a bad vertex is never added.
    Since $E_{j-1}$ is determined by $\cF_{i_j-1}$, it follows that
    \[
        \pr[E_j]
        =\E\big[\mathbf{1}_{E_{j-1}}\pr[v_{i_j}\in I\mid\cF_{i_j-1}]\big]
        \le\hp\,\pr[E_{j-1}].
    \]
    Applying this inequality for $j=1,\ldots,s$ gives $\pr[S\subseteq I] = \pr[E_s] \le\hp^s$.
\end{proof}

\subsection{The probability of containing a given vertex}

We now turn to property~\ref{item: vertex prob}.
By Lemma~\ref{lemma: K3 martingale}, the expected weight of a vertex just before it is processed is $\phi$, so the only way for a vertex to be added with probability less than $\phi$ is for it to become bad.
To show that this is unlikely, we track the \emph{entropy contribution} at $v_j$ after the $i$-th step, defined by
\[
    h_i(v_j) \coloneqq -p_i(v_j)\log p_i(v_j),
\]
with the usual convention that $0 \log 0 = 0$.
The following lemma bounds the expected entropy contribution at $v_j$ just before $v_j$ is processed.
We defer its proof to Section~\ref{subsec: entropy loss} and first show how it implies property~\ref{item: vertex prob}.

\begin{lemma}\label{lemma: entropy}
    For every $j \in [m]$, we have $\E[h_{j-1}(v_j)] \ge h_0(v_j) - (1+\hp)\Delta\phi^2$.
\end{lemma}

A vertex $v_j$ is added to $I$ precisely when $v_j \notin B_{j-1}$ and $a_j = 1$.
Therefore, by Lemma~\ref{lemma: K3 martingale},
\begin{equation}\label{eq: lb on prob}
    \pr[v_j\in I] = \E\big[p_{j-1}(v_j)\,\mathbf{1}\set{v_j \notin B_{j-1}}\big]
    = \phi - \E\big[p_{j-1}(v_j)\,\mathbf{1}\set{v_j \in B_{j-1}}\big].
\end{equation}
We bound the second term as follows:
\begin{align*}
    \E\big[p_{j-1}(v_j)\,\mathbf{1}\set{v_j \in B_{j-1}}\big]
    &\le \frac{\E\big[ p_{j-1}(v_j)\, \mathbf{1}\set{v_j \in B_{j-1}} \log (p_{j-1}(v_j)/\phi)\big]}{\log (\hp/\phi)} \\
    &\le \frac{\E\big[p_{j-1}(v_j)\log (p_{j-1}(v_j)/\phi)\big]}{\log(\hp/\phi)} \\
    &= \frac{\E\big[p_{j-1}(v_j)\big]\log(1/\phi) + \E\big[p_{j-1}(v_j)\log p_{j-1}(v_j)\big]}{\log(\hp/\phi)} \\
    &= \frac{h_0(v_j) - \E\big[h_{j-1}(v_j)\big]}{\log(\hp/\phi)}
    \le \frac{(1+\hp)\Delta\phi^2}{\log(\hp/\phi)}.
\end{align*}
The first inequality holds because $p_{j-1}(v_j) \ge \hp$ whenever $v_j \in B_{j-1}$.
The second holds because every weight is either $0$ or at least $\phi$, so that $p_{j-1}(v_j)\log(p_{j-1}(v_j)/\phi)\ge0$ (interpreting the left-hand side as $0$ when the weight vanishes).
The second equality uses $\E[p_{j-1}(v_j)] = \phi$ (by Lemma~\ref{lemma: K3 martingale}) and $h_0(v_j) = \phi\log(1/\phi)$, and the final inequality is Lemma~\ref{lemma: entropy}.
Combining this bound with~\eqref{eq: lb on prob} yields
\[
    \pr[v_j \in I] \ge \phi - \frac{(1+\hp)\Delta\phi^2}{\log(\hp/\phi)} = \phi\left(1 - \frac{(1+\hp)\Delta\phi}{\log(\hp/\phi)}\right),
\]
which establishes property~\ref{item: vertex prob}.
Together with Lemmas~\ref{lemma: 2 ISET prob} and~\ref{lemma: all sets prob}, this completes the proof of Proposition~\ref{prop: triangle free}, modulo Lemma~\ref{lemma: entropy}.

\subsection{Controlling the entropy loss}\label{subsec: entropy loss}

It remains to prove Lemma~\ref{lemma: entropy}.
To this end, for distinct vertices $u,v\in V(G)$, we define the \emph{weight product} after the $i$-th step as
\[
    \xi_i(uv) \coloneqq p_i(u)p_i(v).
\]
Intuitively, processing a vertex $v_k \in N_L(v_j)$ decreases the expected entropy contribution at $v_j$ by at most a constant multiple of the current weight product $\xi_i(v_kv_j)$.
The following lemma makes this precise.

\begin{lemma}\label{lemma: entropy wrt i}
    For all $j \in [m]$ and $0 \leq i \leq j - 1$, we have
    \[
        \E[h_{j-1}(v_j) \mid \cF_i] \geq h_i(v_j) - (1+\hp)\sum_{\substack{v_k \in N_L(v_j) \\ k > i}}\xi_i(v_kv_j).
    \]
\end{lemma}

\begin{proof}
\stepcounter{ForClaims}
\renewcommand{\theForClaims}{\ref*{lemma: entropy wrt i}}
    We proceed by backward induction on $i$.
    The case $i = j-1$ is trivial, since the sum is empty and both sides equal $h_{j-1}(v_j)$.
    Now suppose that the claim holds for some $1 \leq i \leq j-1$; we prove it for $i - 1$.
    By the tower property and the induction hypothesis,
    \begin{align}
        \E[h_{j-1}(v_j) \mid \cF_{i-1}] &= \E\big[\E[h_{j-1}(v_j) \mid \cF_i] \bigm| \cF_{i-1}\big] \nonumber \\
        &\ge \E\Bigg[ h_i(v_j) - (1+\hp)\sum_{\substack{v_k \in N_L(v_j) \\ k > i}}\xi_i(v_kv_j) \,\Bigg|\, \cF_{i-1}\Bigg] \nonumber \\
        &= \E[ h_i(v_j)\mid \cF_{i-1}] - (1+\hp)\sum_{\substack{v_k \in N_L(v_j) \\ k > i}}\E[\xi_i(v_kv_j)\mid \cF_{i-1}].\label{eq: induction}
    \end{align}
    We bound the two terms separately.

    \begin{claim}\label{claim: entropy expectation}
        $\E[ h_i(v_j)\mid \cF_{i-1}] \ge h_{i-1}(v_j) - (1+\hp)\,\mathbf{1}\set{v_i \in N_L(v_j)}\,\xi_{i-1}(v_iv_j)$.
    \end{claim}
    \begin{claimproof}
        If $v_i \notin N_L(v_j)$, $v_i \in B_{i-1}$, or $v_j \in B_{i-1}$, then $h_i(v_j) = h_{i-1}(v_j)$ and the claim is trivial.
        We may therefore assume that $v_i \in N_L(v_j)$ and $\set{v_i, v_j} \cap B_{i - 1} = \emptyset$.
        Then
        \begin{align*}
            \E[ h_i(v_j)\mid \cF_{i-1}] &= -p_{i - 1}(v_i)\cdot (0\log 0) - (1-p_{i-1}(v_i))\cdot \frac{p_{i-1}(v_j)}{1 - p_{i-1}(v_i)}\log\left(\frac{p_{i-1}(v_j)}{1 - p_{i-1}(v_i)}\right) \\
            &= h_{i-1}(v_j) - p_{i-1}(v_j)\log\left(\frac{1}{1 - p_{i-1}(v_i)}\right).
        \end{align*}
        Since $v_i \notin B_{i-1}$, we have $a \coloneqq p_{i - 1}(v_i) < \hp \le 1/2$, and so
        \[
            -\log(1-a)\le a+\frac{a^2}{2(1-a)}\le a+a^2\le(1+\hp)a.
        \]
        Therefore,
        \[
            \E[ h_i(v_j)\mid \cF_{i-1}] \ge h_{i-1}(v_j) - (1+\hp)p_{i - 1}(v_i)p_{i-1}(v_j) = h_{i-1}(v_j) - (1+\hp)\,\xi_{i-1}(v_iv_j),
        \]
        as desired.
    \end{claimproof}

    \begin{claim}\label{claim: energy expectation}
        For every $v_k \in N_L(v_j)$ with $k > i$, we have $\E[\xi_i(v_kv_j)\mid \cF_{i-1}] \le \xi_{i-1}(v_kv_j)$.
    \end{claim}
    \begin{claimproof}
        Since $v_kv_j \in E(G)$ and $G$ is triangle-free, we have $N_L(v_j) \cap N_L(v_k) = \emptyset$.
        Thus at most one of the weights of $v_j$ and $v_k$ changes in the $i$-th step.
        If neither changes, the claim is trivial.
        If one of them changes, the other is determined by $\cF_{i-1}$, and the argument of Lemma~\ref{lemma: K3 martingale} applied to the changing weight gives the claim with equality.
    \end{claimproof}

    Plugging Claims~\ref{claim: entropy expectation} and~\ref{claim: energy expectation} into~\eqref{eq: induction} gives
    \[
        \E[h_{j-1}(v_j) \mid \cF_{i-1}] \ge h_{i-1}(v_j) - (1+\hp)\sum_{\substack{v_k \in N_L(v_j) \\ k > i-1}}\xi_{i-1}(v_kv_j),
    \]
    which completes the induction.
\end{proof}

\begin{proof}[Proof of Lemma~\ref{lemma: entropy}]
    Apply Lemma~\ref{lemma: entropy wrt i} with $i = 0$.
    Since $\cF_0$ is trivial and $\xi_0(v_kv_j) = \phi^2$ for every edge $v_kv_j$, we obtain
    \[
        h_0(v_j) - \E\big[h_{j-1}(v_j)\big] \le (1+\hp)\sum_{v_k \in N_L(v_j)}\phi^{2} = (1+\hp)\,d_L(v_j)\,\phi^{2} \le (1+\hp)\Delta\phi^{2},
    \]
    as claimed.
\end{proof}

\subsection*{Acknowledgements}

We thank Igor Araujo, Ce Chen, Patryk Morawski, Jacques Verstra\"{e}te, Minh-Quan Vo, and Chase Wilson for insightful discussions, particularly about adapting the Shearer-type approach of Verstra\"{e}te and Wilson~\cite{VW} for uncrowded hypergraphs.
We are also grateful to Jacques Verstra\"{e}te for generously sharing his ideas during his visit to the University of Illinois Urbana--Champaign in April 2025.

\subsection*{Statement on the use of AI}

The authors used GPT Astra to assist with computations involving triangles in the random independent set $I$ sampled by Algorithm~\ref{algorithm: fcp k3}, including identifying suitable $2$-independent configurations and using the algorithm's probability bounds to estimate how many survive. The authors then rewrote and carefully revised these parts. GPT Astra and Claude Opus 5.5 were also used to polish parts of the paper and identify typos. 
The authors developed all other mathematical ideas, including the key idea of sampling an independent set $I$ in the graph of \emph{heavy} pairs using Algorithm~\ref{algorithm: fcp k3} and then deleting vertices from triangles contained in $I$ to obtain a large triangle-free induced subhypergraph to which the Li--Postle theorem (Theorem~\ref{thm:LP}) applies. The authors take full responsibility for the content of the paper.

\bibliographystyle{amsplain}
\bibliography{KMV_proof}

\appendix
\numberwithin{equation}{section}
\section{Proof of Lemma~\ref{lem:compression}}\label{app:compression}

For the reader's convenience, we restate the lemma.

\begin{lemma*}[Restatement of Lemma~\ref{lem:compression}]
    For every integer $r\ge3$, there is a constant $C_r>0$ such that the following holds. If $\cH$ is an $n$-vertex hypergraph with no loose triangle in which every edge has size between $3$ and $r$, then there is a hypergraph $\cH'$ on $V(\cH)$ with the following properties:
    \begin{enumerate}[label=\textup{(\roman*)}]
        \item $\cH'$ has no loose triangle, has edge sizes between $3$ and $r$, and is an antichain;
        \item every independent set of $\cH'$ is independent in $\cH$;
        \item every triple $X\subseteq V(\cH')$ satisfies $\dd_4(X,\cH')\le 2r$; and
        \item $m_j(\cH')\le C_rn^2$ for every $4\le j\le r$.
    \end{enumerate}
\end{lemma*}

The argument follows~\cite[Lemmas~3.2--3.4, Corollary~3.5, and Lemma~3.6]{KMV}.
All hypergraphs below have vertex set $V(\cH)$ and edge sizes between $3$ and $r$.
Recall from Definition~\ref{def: expandable} that a nonempty set $X$ is \emph{expandable} in a hypergraph $\cK$ if, for every $Y\subseteq V(\cH)\setminus X$ with $|Y|\le2r$, some edge of $\cK$ contains $X$ and is disjoint from $Y$.
Every edge of $\cK$ is expandable in $\cK$, and every expandable set lies in an edge and hence has size at most $r$.
By Lemma~\ref{lem:expand}, if $\cK$ contains no loose triangle, then no three sets expandable in $\cK$ form a loose triangle.
\begin{proof}[Proof of Lemma~\ref{lem:compression}]
\stepcounter{ForClaims}
\renewcommand{\theForClaims}{\ref*{lem:compression}}
Starting from $\cK=\cH$, we repeatedly replace an edge $e$ by a proper subset $X\subsetneq e$ of size at least three that is expandable in the current hypergraph $\cK$.
Each replacement produces a hypergraph $\cK^*$ on the same vertex set, with
\[
    E(\cK^*)\coloneqq\big(E(\cK)\setminus\{e\}\big)\cup\{X\}.
\]
We continue until no such pair $X,e$ exists.
Each replacement strictly decreases the nonnegative integer $\sum_{f\in E(\cK)}|f|$: either $e$ is replaced by the smaller edge $X$, or $X$ is already an edge and $e$ is simply removed.
Thus only finitely many replacements are possible.
Let $\cH'$ be the hypergraph obtained when no further replacement is possible.

\begin{claim}\label{claim: appendix structure}
    The hypergraph $\cH'$ satisfies \ref{item: compression structure} and~\ref{item: compression independent}.
    Moreover, if a set $X$ of at least three vertices is properly contained in an edge of $\cH'$, then $X$ is not expandable in $\cH'$.
\end{claim}

\begin{claimproof}
    Consider a replacement of an edge $e\in E(\cK)$ by a set $X$ of at least three vertices which is properly contained in $e$ and expandable in $\cK$, producing $\cK^*$ as above.
    Every edge of $\cK^*$ is expandable in $\cK$: each edge other than $X$ is an edge of $\cK$, and $X$ is expandable by its choice.
    Hence, by Lemma~\ref{lem:expand}, if $\cK$ has no loose triangle, then neither does $\cK^*$.
    The edge sizes remain between $3$ and $r$, since $3\le|X|<|e|\le r$.
    An independent set of $\cK^*$ contains no edge in $E(\cK)\setminus\{e\}$ and does not contain $X$.
    Since $X\subseteq e$, it cannot contain $e$ either, so it is independent in $\cK$.
    Applying these observations after each replacement shows that $\cH'$ has no loose triangle, has edge sizes between $3$ and $r$, and satisfies~\ref{item: compression independent}.

    If a set $X$ of size at least three were expandable in $\cH'$ and properly contained in an edge $e$ of $\cH'$, we could replace $e$ by $X$ once more, contrary to the choice of $\cH'$ as a hypergraph in which no further replacement is possible.
    In particular, an edge of $\cH'$ cannot be properly contained in another edge, since it has size at least three and is expandable in $\cH'$.
    Thus $\cH'$ is an antichain.
\end{claimproof}

\begin{claim}\label{claim: appendix codegrees}
    For all integers $3\le i\le j\le r$ and every set $X\subseteq V(\cH')$ of size $i$, we have $\dd_j(X,\cH')\le(2r)^{j-i}$.
\end{claim}

\begin{claimproof}
    Fix $j\in\{3,\ldots,r\}$.
    We prove the bound by downward induction on $i$, starting at $i=j$ and ending at $i=3$.
    When $i=j$, the only possible edge of size $j$ containing $X$ is $X$ itself, so $\dd_j(X,\cH')\le1=(2r)^{j-i}$.

    Let $3\le i<j$, and suppose that the bound holds for every set of size $i+1$.
    Fix a set $X\subseteq V(\cH')$ of size $i$.
    If $\dd_j(X,\cH')=0$, the bound is immediate.
    Otherwise, $X$ is properly contained in an edge of size $j$, so it is not expandable by Claim~\ref{claim: appendix structure}.
    Hence some set $Y\subseteq V(\cH')\setminus X$ with $|Y|\le2r$ meets every edge containing $X$.
    Every edge of size $j$ containing $X$ therefore contains $X\cup\{y\}$ for some $y\in Y$.
    Since $|X\cup\{y\}|=i+1$ for every $y\in Y$, the induction hypothesis gives $\dd_j(X\cup\{y\},\cH')\le(2r)^{j-i-1}$.
    Therefore,
    \[
        \dd_j(X,\cH')
        \le\sum_{y\in Y}\dd_j(X\cup\{y\},\cH')
        \le |Y|(2r)^{j-i-1}
        \le(2r)^{j-i},
    \]
    as claimed.
\end{claimproof}

Taking $i=3$ and $j=4$ in Claim~\ref{claim: appendix codegrees} proves \ref{item: compression codegree} when $r\ge4$; when $r=3$, there are no edges of size four, so $\dd_4(X,\cH')=0$ for every triple $X$ and the desired bound holds trivially.

\begin{claim}\label{claim: appendix pairs}
    If a pair $P\subseteq V(\cH')$ is not expandable in $\cH'$, then $\dd_j(P,\cH')\le(2r)^{j-2}$ for every $3\le j\le r$.
\end{claim}

\begin{claimproof}
    Fix $j\in\{3,\ldots,r\}$.
    Since $P$ is not expandable, some set $Y\subseteq V(\cH')\setminus P$ with $|Y|\le2r$ meets every edge containing $P$.
    Thus every edge of size $j$ containing $P$ contains a triple $P\cup\{y\}$ for some $y\in Y$.
    Applying Claim~\ref{claim: appendix codegrees} with $i=3$ to these triples gives
    \[
        \dd_j(P,\cH')
        \le\sum_{y\in Y}\dd_j(P\cup\{y\},\cH')
        \le |Y|(2r)^{j-3}
        \le(2r)^{j-2},
    \]
    as claimed.
\end{claimproof}

\begin{claim}\label{claim: appendix pair in edge}
    Every edge of $\cH'$ contains a pair that is not expandable in $\cH'$.
\end{claim}

\begin{claimproof}
    Fix an edge $e\in E(\cH')$ and choose distinct vertices $x,y,z\in e$, which is possible since $|e|\ge3$.
    If the pairs $\{x,y\}$, $\{x,z\}$, and $\{y,z\}$ were all expandable in $\cH'$, they would form a loose triangle of expandable sets.
    This contradicts Lemma~\ref{lem:expand}, since $\cH'$ has no loose triangle.
    Hence at least one of these pairs is not expandable in $\cH'$.
\end{claimproof}

Let $\cN$ be the set of pairs of vertices that are not expandable in $\cH'$.
By Claim~\ref{claim: appendix pair in edge}, every edge of $\cH'$ contains at least one pair in $\cN$.
Thus, for every $4\le j\le r$, counting edges through these pairs and applying Claim~\ref{claim: appendix pairs} gives
\[
    m_j(\cH')
    \le\sum_{P\in\cN}\dd_j(P,\cH')
    \le |\cN|(2r)^{j-2}
    \le\binom n2(2r)^{j-2}
    \le\frac{(2r)^{r-2}}2\,n^2.
\]
Taking $C_r\coloneqq(2r)^{r-2}/2$ proves \ref{item: compression sparse} and completes the proof.
\end{proof}

\end{document}